\documentclass[10pt]{article}
\usepackage{graphicx}
\usepackage{biograph}        
\usepackage{amsmath,amsfonts,amssymb,amsthm}
\usepackage{latexsym,xypic,eucal,mathrsfs,setspace}
\newcommand{\n}{\mathfrak{N}}
\newcommand{\g}{\mathfrak{g}}

\newcommand{\z}{\mathfrak{Z}}
\newcommand{\w}{\mathfrak{W}}
\newcommand{\cart}{\mathfrak{h}}
\newcommand{\chev}{\mathcal{C}}
\newcommand{\chevo}{\mathcal{C}_{{0}}}
\newcommand{\goo}{\mathfrak{g}_{{0}}}
\newcommand{\goc}{\g_{{0}}^{{\mathbb{C}}}}
\newcommand{\ngam}{\Gamma \backslash N}
\newcommand{\Rq}{\mathbb{R}^q}
\newcommand{\so}{\mathfrak{so}(q,\mathbb{R})}
\newcommand{\ip}{\langle \ , \ \rangle}
\newcommand{\vv}{\mathfrak{V}}
\newcommand{\carts}{\mathfrak{h}^*}
\newcommand{\cartr}{\mathfrak{h}_{{\mathbb{R}}}}
\newcommand{\weyl}{\mathcal{W}}
\newcommand{\linl}{\lambda\in\Lambda}
\newcommand{\aone}{\alpha_1}
\newcommand{\atwo}{\alpha_2}
\newcommand{\ai}{\alpha_i}
\newcommand{\an}{\alpha_n}

\newcommand{\afour}{\alpha_4}

\newcommand{\mone}{m_1}
\newcommand{\mtwo}{m_2}
\newcommand{\mn}{m_n}

\newcommand{\mthree}{m_3}
\newcommand{\mfour}{m_4}

\newcommand{\mi}{m_i}

\newcommand{\curvegreat}{\succ}
\newcommand{\carto}{\mathfrak{h}_0}
\newcommand{\cartoc}{\mathfrak{h}_{{0}}^{{\mathbb{C}}}}
\newcommand{\klam}{K_{\lambda,\mu}}
\newcommand{\kzeta}{K_{\zeta,\mu}}
\newcommand{\ej}{e^{tj(\alpha Z)}}
\newcommand{\ejo}{e^{\omega j(\alpha Z)}}
\newcommand{\jinv}{j(\alpha Z)^{-1}}
\newcommand{\jjinv}{j( Z)^{-1}}
\newcommand{\imdf}{\mbox{Im}(dF_Z)}

\newtheorem{theorem}{Theorem}[section]
\newtheorem{lemma}[theorem]{Lemma}
\newtheorem{corollary}[theorem]{Corollary}
\newtheorem{proposition}[theorem]{Proposition}

\newtheorem{definition}[theorem]{Definition}
\newtheorem{sublemma}{Lemma}[subsection]

\newtheorem{example}[theorem]{Example}
\newtheorem{remark}[theorem]{Remark}

\begin{document}

\title{Closed Geodesics on Compact Nilmanifolds with Chevalley Rational Structure
}


\author{Rachelle C. DeCoste    \\
             Department of Mathematical Sciences\\ The United States Military Academy\\ West Point, NY 10996 \\
              Tel.: 845-938-2530\\
              Rachelle.DeCoste@usma.edu   
}



\date{\today}

\maketitle

\begin{abstract}
We continue the study of the distribution of closed geodesics on nilmanifolds $\ngam$, constructed from a simply connected 2-step nilpotent Lie group $N$ with a left invariant metric and a lattice $\Gamma$ in $N$. We consider a Lie group $N$ with associated  2-step nilpotent Lie algebra $\n$ constructed from an irreducible representation of a compact semisimple Lie algebra $\goo$ on a real finite dimensional vector space $U$.  We determine sufficient conditions on the semisimple Lie algebra $\goo$ for $\ngam$ to have the density of closed geodesics property where $\Gamma$ is a lattice arising from a Chevalley rational structure on $\n$.
\end{abstract}

\section{Introduction}
\label{intro}
In this paper we continue the study of the distribution of closed geodesics on compact 2-step nilmanifolds.  This study was begun by Eberlein \cite{eb1}, Mast \cite{mast} and Lee-Park \cite{leepark} and continued more recently by DeMeyer \cite{demeyer}.  

The objects we study are analogous to tori. Considering the torus as $T^n=\mathbb{R}^n/\mathbb{Z}^n$, the smoothly closed geodesics on $T^n$ are easily described as projections of  those straight lines in $\mathbb{R}^n$ that pass through points of the lattice $\mathbb{Z}^n$.  It is straightforward to see that  these closed geodesics are dense in all the geodesics on $T^n$.  In general, a compact manifold is said to have the density of closed geodesics property (DCG) if the vectors tangent to the  unit speed geodesics are dense in the unit tangent bundle.  We consider nilmanifolds $\ngam$ constructed from a 2-step nilpotent Lie group $N$ with a lattice $\Gamma$.  Our method of proving the DCG for manifolds $\ngam$ depends on showing the density of curves passing through lattice points.
We construct $N$ from an irreducible representation of a compact semisimple Lie algebra $\goo$ on a real finite dimensional vector space $U$.  The left invariant metric on $N$ induces a metric on $\ngam$.

Lee and Park have shown that if $\ngam$ satisfies a nonsingularity condition and a resonance condition then the density of closed geodesics property will always hold.  DeMeyer first investigated whether the nonsingularity condition was necessary in the case where $\goo=\mathfrak{su}(2)$.  She considered the irreducible, real representations of $\mathfrak{su}(2)$ on an odd dimensional space.  We investigate a general case $\n=U\oplus\goo$ where the nonsingularity condition does not hold for a compact semisimple Lie algebra $\goo$ and some $\goo-$module $U$.   We determine sufficient conditions on the associated semisimple Lie algebra $\goo$ for $\ngam$ to have the density of closed geodesics property for $\Gamma$ determined by a compact Chevalley basis of $\goo$.  
 We expand upon previous results in the following main theorem, where  the $\goo-$module $U$ is defined to be admissible if the roots of $\g=\goc$ are weights of $V=U^\mathbb{C}$ that satisfy a certain multiplicity condition as discussed in Section \ref{sec3.3}.  We show that almost all $\goo-$modules $U$ are admissible and we obtain many new nilmanifolds that satisfy the DCG property.

\begin{theorem}[Main Result]  
Let $U$ be an admissible $\g_0-$module for a compact semisimple  Lie algebra $\goo$.  Let $\n=U\oplus\goo$ be the metric 2-step nilpotent Lie algebra defined in Section \ref{chap4} and let $N$ be the corresponding simply connected, metric 2-step nilpotent Lie group.  
Let $\beta$ be a rational basis of $\n$ determined by a Chevalley basis of $\g=\g_0^\mathbb{C}$.
Then $\Gamma\backslash N$ satisfies DCG for every lattice $\Gamma$ in $N$ determined by $\beta$.
\end{theorem}

\subsection{Definitions and History}
\label{sec1.1}

Let $\n$ be a Lie algebra.  In our examination of 2-step nilpotent Lie algebras we are considering
those nonabelian Lie algebras which are almost abelian; i.e. $[\n,\n]\subseteq \z$, the center.
We call a simply connected Lie group $N$ 2-step nilpotent if its associated Lie algebra $\n$ is 2-step nilpotent.    Denote by $\exp: \n\rightarrow N$ the Lie group exponential map with inverse $\log: N \rightarrow \n$.  It is known that the map $\exp: \n\rightarrow N$  is a diffeomorphism (cf. \cite{rag}).

We define a 2-step nilpotent Lie algebra $\n$ to be of type $(p,q)$ if $\dim [\n,\n]=p$
and $\mbox{codim}[\n,\n]=q$.  Every 2-step nilpotent Lie algebra is isomorphic to
one of the following metric examples.  Let $\w$ be a $p-$dimensional subspace of $\so$, the 
$q\times q$ real symmetric matrices, and let $\n=\Rq\oplus\w$.  We endow $\n$ with the left
invariant metric such that $\Rq$ and $\w$ are orthogonal, $\Rq$ has the standard inner product
and $\langle Z_1, Z_2\rangle=-\mbox{trace}(Z_1Z_2)$ for all $Z_1, Z_2\in\w$.  Then 
the bracket relations on $\n$ are defined by letting $\w$ be 
contained in the center of $\n$, $[\Rq,\Rq]\subseteq\w$ and by requiring that
$\langle [X,Y],Z\rangle=\langle Z(X),Y\rangle$ for all $X, \ Y$ in $\Rq$.  Interesting cases arise when $\w$ is a subalgebra of $\so$ 
or  a Lie triple system:  $[\w,[\w,\w]]\subseteq\w$.  We focus on the former letting $\w=\goo$ a compact, semisimple Lie algebra with $\goo-$module $U$ with inner product $\ip$ for which the elements of $\goo$ are skew symmetric.

\subsection{Lattices and Rational Structures}
\label{sec1.2}

Let $N$ be a simply connected 2-step nilpotent Lie group with a left invariant metric.  A discrete subgroup $\Gamma$ of $N$ is a lattice
in $N$ if $\Gamma \backslash N$ admits a finite invariant measure.  Equivalently by \cite[Theorem 2.1]{rag}, a discrete subgroup $\Gamma$ of a simply connected nilpotent Lie group $N$ is a lattice if and only if $\ngam$ is compact.

Not every Lie algebra admits a lattice; see \cite{rag} section 2.14 for the construction of an example.  Mal'cev described a criterion on the associated Lie algebra for a Lie group $N$ to admit a lattice.

\begin{theorem}\cite{malcev} Let $N$ be a simply connected, nilpotent Lie group, and let $\n$ be its Lie algebra.  Then $N$ admits a lattice $\Gamma$ if and only if $\n$ admits a basis $\beta=\{X_1,\dots,X_n\}$ such that $$\displaystyle [X_i,X_j]=\sum_{k=1}^{n}C_{ij}^k X_k$$ for all $i,j$ where the constants $\{C_{ij}^k\}$ are all rational.  \label{malcev}\end{theorem}

The set $\{C_{ij}^k\}$ is called the set of structure constants of the basis $\beta$.  The basis $\beta$ is called a rational basis if the structure constants are all rational.  In this case we say that $\n_\mathbb{Q}=\mathbb{Q}-\mbox{span}\{\beta\}$ is a rational structure for $\n$.  
We define a  vector $X$ to be rational if $X\in \mathbb{Q}-\mbox{span}(\beta)$ where $\beta$ is a rational basis.  A subspace $\w\subseteq \n$ is rational
if $\w$ has a basis of rational vectors. If a subalgebra $\w$ is rational as a subspace, then we say that the group $W=\exp(\w)$ is a rational subgroup of $N$.

The following results establish the relationship between lattices in a Lie group $N$ and rational structures for the associated Lie algebra $\n$.  

\begin{theorem}\cite{cg}  Let $N$ be a nilpotent Lie group with Lie algebra $\n$.  
\begin{enumerate}
\item If $N$ admits a lattice $\Gamma$, then $\n_{\mathbb{Q}}=\mathbb{Q}-\mbox{span}\{\log\Gamma\}$ is a rational structure for $\n$.
\item If $\beta$ is a basis for $\n$ with rational structure constants, then the group $\Gamma$ in $N$ generated by $\exp(\mathbb{Z}-\mbox{span}(\beta))$ is a lattice in $N$.  Moreover, $\mathbb{Q}-\mbox{span}\{\log\Gamma\}=\mathbb{Q}-\mbox{span}(\beta)$.
\end{enumerate}
\end{theorem}

\begin{proposition}\cite[Proposition 5.3]{eb1}, Let $\mathcal{Z}$ be the center of $N$ and let $\Gamma$ be a lattice in $N$. 
 Let the map $\log: N \rightarrow \n$ denote the inverse
of the exponential map $\exp: \n\rightarrow N$.  Then $\Gamma \cap \mathcal{Z}$ is a lattice subgroup in $\mathcal{Z}$ and $(\log\Gamma)\cap\z$ is a vector lattice in $\z$.\label{eb15.3}
\end{proposition}

Recall that a Lie subgroup $H$ of $N$ is rational if its Lie algebra $\cart$ is a rational subalgebra of $\n$.

\begin{proposition}\cite[Theorem 5.1.11]{cg} Let $\Gamma$ be a lattice in $N$.  A subgroup $H$ of $N$ is rational with respect to $\n_\mathbb{Q}=\mathbb{Q}-\mbox{span}\{\log \Gamma\}$ if and only if $H\cap \Gamma$ is a lattice in $H$.
\label{cglattice}\end{proposition}

\begin{example}[Chevalley Rational Structure]
As above we define a metric, 2-step nilpotent Lie algebra $\n=U\oplus\goo$, where $\goo$ is a compact, semisimple Lie algebra and $U$ is a finite dimensional real $\goo-$module.  A Chevalley basis $\chev$ for $\g=\goc$ defines a compact Chevalley basis $\chevo$ for $\goo$.  A result of Raghunathan \cite{ragart} states that $U$ admits a basis $\mathcal{B}$ such that $\chevo$ leaves invariant $\mathbb{Q}-\mbox{span}(\mathcal{B})$.  If the $\goo-$invariant inner product is chosen suitably on $\n$, then $\n_\mathbb{Q}=\mathbb{Q}-\mbox{span}\{\mathcal{B}\cup\chevo\}$ is a rational structure on $\n$ that we call a Chevalley rational structure.
\end{example}

The following concept is used in the proof of our main result.  We say that two lattices $\Gamma_1$ and $\Gamma_2$ are commensurable if $\Gamma_1\cap\Gamma_2$ has a finite index in $\Gamma_1$ and $\Gamma_2$.  
\begin{remark}
\label{remark1.3.4}
By Theorem 5.4.2 of \cite{cg}, for any lattice $\Gamma$ there exist lattices $\Gamma_1$ and $\Gamma_2$ commensurable with $\Gamma$ such that $\log \Gamma_1$ and $\log \Gamma_2$ are vector lattices in $\n$ and $\Gamma_1\subseteq\Gamma\subseteq\Gamma_2$.  Hence if $r>0$ is sufficiently large, then any ball in $\n$ of radius $r$ intersects $\log \Gamma$.
\end{remark}
\subsection{Density of Closed Geodesics}
\label{sec1.3}


Let $\n$ be a metric 2-step nilpotent Lie algebra and write $\n=\vv\oplus\z$, where $\z$ is the center of $\n$ and $\vv=\z^\perp$.  
For every
nonzero element $Z\in\z$, we define a skew symmetric linear transformation $j:\vv\rightarrow \vv$ 
which satisfies
$$\langle [X,Y],Z\rangle=\langle j(Z)X,Y\rangle$$ for all $X, \ Y$ in $\vv$.  

For the following, fix an inner product $\langle \ , \ \rangle$ on $\n$.  

\begin{definition}
The Lie algebra $\{\n,\langle \ , \ \rangle\}$ is \emph{nonsingular} if $j(Z)$ is a nonsingular map for all nonzero $Z\in\z$.  The Lie algebra $\{\n,\langle \ , \ \rangle \}$ is \emph{almost nonsingular} if $j(Z)$ is nonsingular for all $Z$ in some open dense subset of $\z$.  Lastly, $\{\n, \langle \ , \ \rangle\}$ is \emph{singular} if $j(Z)$ has a nonzero kernel for all $Z\in\z$.
\end{definition}

The condition of $\n$ being almost nonsingular is equivalent to $j(Z)$ being nonsingular for some $Z\in\z$ (cf. \cite{leepark}).  
Nonsingularity and almost nonsingularity are conditions independent of the metric.  For example, the Lie algebra $\{\n, \langle \ , \ \rangle\}$ is nonsingular if and only if $\mbox{ad}(X):\n \rightarrow \z$ is surjective for every $X\in\n-\z$.  This and the observation above from \cite{leepark} lead to the following result.

\begin{proposition}\cite{gm} Every 2-step nilpotent Lie algebra is either nonsingular, almost nonsingular, or singular.
\end{proposition}

\paragraph{Geodesics}
Let $N$ be a metric 2-step nilpotent Lie group with associated Lie algebra $\n=\vv\oplus\z$.  To describe the geodesics of $N$, it suffices to consider those geodesics that begin at the identity of $N$ since $N$ has a left invariant metric.  Let $\gamma(t)$ be a curve in $N$ with $\gamma(0)=e$ and $\gamma'(0)=X_0+Z_0\in\n$ for some $X_0\in \vv$ and $Z_0\in\z$.  Then we write $\gamma(t)=\exp (X(t)+Z(t))$ with $X(t)\in \vv$ and $Z(t)\in\z$ for all $t$ and $X'(0)=X_0$, $Z'(0)=Z_0$.  The functions $X(t), \ Z(t)$ are unique since $\exp: \n\rightarrow N$ is a diffeomorphism.

A. Kaplan \cite{kaplan} shows that a curve $\gamma(t)$ is a geodesic in $N$ if and only if the following equations are satisfied:

\begin{enumerate}
\item $X''(t)=j(Z_0)X'(t)$,
\item $Z'(t)+\frac{1}{2}[X'(t),X(t)]\equiv Z_0$ for all $t\in\mathbb{R}$.
\end{enumerate}

\begin{definition} A compact Riemannian manifold $M$ is said to satisfy the \emph{density of closed geodesics property} (DCG property)
if the vectors tangent to unit speed geodesics are dense in the unit tangent bundle of $M$. \end{definition}

For $\n=\vv\oplus\z$, if $Z$ is a nonzero element of $\z$, then the map  $j(Z):\vv\rightarrow \vv$ is defined as above for $X,Y \in \vv$.  Since $j(Z)$ is a skew symmetric map, its eigenvalues are purely imaginary.  Thus the ratio of any two eigenvalues is real.  
\begin{definition} We say that the map $j(Z)$ is in 
\emph{resonance} if the ratio of any two nonzero eigenvalues of $j(Z)$ is a rational real number. 
\end{definition}
 Eberlein \cite{eb1} describes another characterization of resonance in the following result.

\begin{lemma}\cite[Lemma 4.23]{eb1}, Let $Z$ be any nonzero element of $\z$.  Then $j(Z)$ is in resonance if and only if
$e^{\omega j(Z)}$ is the identity on $\vv$ for some $\omega >0$ where $e$ is the matrix exponential map.\label{onefourtwo}\end{lemma}

As seen in Theorem \ref{mastres} below, Mast \cite{mast} has shown that the density of resonant vectors in $\z$ is necessary for a manifold $\n=\vv\oplus \z$ to have the DCG property.

We show in Proposition \ref{ratlres} that for $\n=U\oplus\goo$ as above, $\goo$ has a dense set of resonant vectors.
It is not difficult to construct examples of nonresonant behavior; see for example,  Eberlein \cite[example 5.8]{eb1}.  

\begin{definition} Let $\phi$ be an arbitrary element of $N$.  We say that $\phi$ \emph{translates} a unit speed geodesic
$\gamma(t)$ in $N$ by an amount $\omega$ if $\phi\cdot \gamma(t)=\gamma(t+\omega)$ for all 
$t\in\mathbb{R}$. \end{definition} 

To prove that a curve $\gamma(t)$ in $N$ projects to a closed geodesic in $\Gamma \backslash N$, it suffices to show that there exists an element of the lattice, $\phi\in\Gamma$, such that $\phi$ translates $\gamma(t)$ in $N$.  Thus, to establish the density of closed geodesics property in $\ngam$, it suffices to show that $N$ contains a dense set of geodesics that are translated by elements of the lattice $\Gamma$.
We relate the property of resonance to the translation of geodesics through the following result of Eberlein.

\begin{proposition} \cite[Proposition 4.3]{eb1}  Let $\gamma(t)$ be a unit speed geodesic with $\gamma(0)=e$ and $\gamma'(0)=X_0+Z_0\in\n=\vv\oplus\z$, where $X_0\in\vv$ and $Z_0\in\z$.  Then the following are equivalent for a positive number $\omega$.
\begin{enumerate}
\item $e^{\omega j(Z_0)}X_0=X_0$.
\item $\gamma(\omega)\cdot \gamma(t)=\gamma(t+\omega)$ for all $t\in\mathbb{R}$.
\end{enumerate}
\label{eb14.3}
\end{proposition}

\subsection{Previous results}
\label{sec1.4}

We include a brief summary of previous results.  All proofs are omitted, but can be found in the original publications as indicated.

Eberlein \cite[Proposition 5.6]{eb1} first proved that for any 2-step nilpotent Lie group $N$ of Heisenberg type  $\ngam$ has the DCG property for any lattice $\Gamma\subseteq N$.  Then for $N$ with a 1-dimensional center, Eberlein found necessary and sufficient conditions for the DCG property to hold.  Lie algebras of  Heisenberg type and Lie algebras with a one-dimensional center are nonsingular Lie algebras.  The following result of Mast concerning the nonsingular case is a generalization of the previous results of Eberlein.

\begin{theorem} \cite{mast} Let $N$ be a nonsingular, simply connected, 2-step nilpotent Lie group with a left invariant metric.
\begin{enumerate} 
\item If $\Gamma\backslash N$ has the density of closed geodesics property for some 
lattice $\Gamma$ then $j(Z)$ is in resonance for a dense subset of $Z\in \z$.
\item If $j(Z)$ is in resonance for all nonzero $Z\in\z$ then for any lattice $\Gamma \subseteq N$,
$\Gamma \backslash N$ has density of closed geodesics.
\end{enumerate}\label{mastres}
\end{theorem}

The proof of this result shows that the resonance condition in (1) is necessary in the singular case as well.

In the following, Lee and Park completely resolve the question of the DCG property for any 2-step nilmanifold associated to a Lie algebra with nonsingular elements, i.e. an almost nonsingular Lie algebra.

\begin{theorem} \cite{leepark}  Let $\n$ be a metric 2-step nilpotent Lie algebra such that
\begin{enumerate}
\item $j(Z)$ is in resonance for a dense set of $Z\in\z$
\item $j(Z)$ is nonsingular for some nonzero $Z\in\z$
\end{enumerate}
Then for any lattice $\Gamma$ of $N$, the nilmanifold $\Gamma\backslash N$ has the density of
closed geodesics property.
\end{theorem}

\paragraph{The singular case}
DeMeyer \cite{demeyer} began the study of Lie algebras in which all $j(Z)$ maps are singular.  It is this case that we continue to consider.

Let $j: \goo\rightarrow \mbox{End}(U)$ be an irreducible real representation of a compact, semisimple Lie algebra $\goo$. Let $\n=U\oplus \goo$ be the corresponding 2-step nilpotent Lie algebra as defined in Section \ref{chap4}.  Here $U$ is equipped with an inner product $\ip$ such that $j(\goo)\subseteq \mathfrak{so}(U,\langle \ , \ \rangle)$.  Note that if $U$ is odd dimensional, then $j(Z)$ has nonzero kernel for every $Z\in\goo$.

Let $V_n$ be the space of complex homogeneous polynomials in two variables of degree $n$.  These are all
of the complex irreducible representations of $\mathfrak{su}(2)$.  If $n$ is even then $V_n=U_n\otimes \mathbb{C}$,  where $U_n$ is a real irreducible $\mathfrak{su}(2)-$module
of real dimension $n+1$.  The following result of DeMeyer concerns these particular real representations.

\begin{theorem}\cite{demeyer} Let $N$ be a simply connected 2-step nilpotent Lie group with a left invariant metric and dimension $2k+4$ that is constructed from an irreducible real representation of $SU(2)$ on a real vector space of dimension $2k+1$, $k\geq 2$.  
Then for any lattice $\Gamma$ in $N$, the compact manifold $\Gamma\backslash N$ with induced metric satisfies the density of closed geodesics property.  \label{demeyer}\end{theorem}

\subsection{Outline of the Proof of our Main Result}
\label{sec1.5}

To prove our main theorem, we rely on the method of proof used by Eberlein, Mast, and DeMeyer which employs a ``first hit map.''  In our case where $\n=U\oplus \goo$ is a real 2-step nilpotent Lie algebra we must consider the complex semisimple Lie algebra $\g=\goc$, the finite dimensional complex $\g-$module $V=U^\mathbb{C}$ and the complex 2-step nilpotent Lie algebra $\n^\mathbb{C}=V\oplus\g$.  Throughout we assume that $\g$ has a nontrivial zero weight space $V_0$.  We find that, with few exceptions, all roots of $\g$ are weights of $V$ with multiplicity greater than or equal to 2.  If our Lie algebra satisfies this root multiplicity condition, then our first hit map $F_Z$, defined for each resonant, rational vector $Z\in\goo$, will have maximal rank on a dense open subset of its domain.  This will allow us to complete the proof of our main result.  We will briefly review basic concepts related to roots and weight spaces of Lie algebras before proving our main result.  Additional review can be found in later sections of this article and in cited sources.

\section{The Lie algebra $\n=U\oplus\goo$}
\label{chap4}

\subsection{The Metric Lie Algebra $\n=U\oplus\goo$}
\label{sec4.1}

Let $\goo$ be a compact, semisimple real Lie algebra, and let $U$ be a finite dimensional real $\goo-$module.  We introduce a 2-step metric nilpotent Lie algebra structure on $\n=U\oplus\goo$, and we study its properties by considering the complex $\g-$module $V=U^\mathbb{C}$, where $\g=\goc$.  See \cite{ebpp} for further details.  Let $\ip$ be an inner product on $U$ such that the elements of $\goo$ are skew symmetric on $U$ relative to $\ip$.  The bracket relations on $\n$ are defined by letting $\goo$ be contained in the center of $\n$, $[U,U]\subseteq \goo$ and by requiring that $\langle [X,Y],Z\rangle=\langle Z(X),Y\rangle$ for all $X,Y\in U$.  Although the $\goo-$invariant inner product $\ip$ on $U$ may vary, the isomorphism type of $\n=U\oplus\goo$ is unchanged (Section 2.4 of \cite{ebpp}).

To study $\n=U\oplus\goo$, we consider the complex semisimple Lie algebra, $\g=\goc$ and the complex $\g-$module $V=U^\mathbb{C}$.  If $\carto$ is a maximal abelian subalgebra of $\goo$, then $\cart=\cartoc$ is the corresponding Cartan subalgebra for $\g$.  

Throughout this paper, let $\g$ denote a complex semisimple Lie algebra and $\mbox{ad}:\g\rightarrow \mbox{End}(\g)$ the adjoint representation. Every semisimple Lie algebra over $\mathbb{C}$ has a Cartan subalgebra $\cart$ and all Cartan subalgebras of $\g$ lie in a single orbit of Aut$(\g)$.

\subsection{Isometries of the Lie Group $N$}
\label{sec4.2}
Let $G_0$ be any compact Lie group, and $\goo$ its Lie algebra.  Let $\n=U\oplus \goo$ be a metric, 2-step nilpotent Lie algebra as described above and let $N$ be the simply connected Lie group with Lie algebra $\n$.   Denote by $I(\n)$ and $I(N)$ the isometry groups of $\n$ and $N$ respectively.  The following result describes the action of $G_0$  on $\n=U\oplus \goo$.

\begin{theorem}\cite[Theorem 3.12]{lauret}  Let $\rho: G_0\rightarrow \mbox{Aut}(U)$ be an irreducible representation of a compact Lie group $G_0$ with discrete kernel on a real, finite dimensional vector space $U$.  Let $\n=U\oplus\goo$ be a metric 2-step nilpotent Lie algebra constructed as above, and let $N$ be the simply connected 2-step nilpotent Lie group with left invariant metric.  Then for each $g\in G_0$ there is a map  $I_g\in \mbox{Aut}(N)\cap I(N)$ such that $dI_g$ acts on $\n=U\oplus \goo$ as an automorphism and an isometry by $(\rho(g),\mbox{Ad}(g))$.\label{laur}\end{theorem}

\section{Roots, Chevalley Basis and Rational Structure}
\label{chap2}

The following review of roots and Chevalley basis can be found in more depth in sources such as Helgason \cite{hel2} and Humphreys \cite{humph}.  We include those definitions and results which are necessary in our further discussion of the DCG of 2-step nilmanifolds.

\subsection{Roots and Root Space Decomposition}
\label{sec2.1}

Let $\Phi\subseteq \carts=\mbox{Hom}(\cart,\mathbb{C})$ denote the set of  all nonzero roots for $\cart$ and let the subset $\Delta$ of $\Phi$  be a base of $\carts$.  Recall that this means that every root $\beta$ can be written as the unique sum $\displaystyle \beta=\sum_{\alpha\in\Delta}m_\alpha \alpha$ where $m_\alpha$ are either all nonnegative or all nonpositive.  The roots in $\Delta$ are called simple roots and if $m_\alpha \geq 0$ for all $\alpha$, then $\beta$ is called a positive root.  Let
 $\Delta\subseteq \Phi$ be a basis of simple roots, $\Delta=\{\alpha_1,\dots,\alpha_r\}$ where 
$r=\dim_\mathbb{C}\cart$ is the rank of $\g$ and of the root system $\Phi$.  Also, let $B$ denote the Killing form of $\g$; i.e. $B(X,Y)=\mbox{Tr}(\mbox{ad }X \mbox{ ad }Y)$. 

\begin{theorem}\cite[Theorem 4.2]{hel2} Let $\Phi$ be the set of nonzero roots.  Then
\begin{enumerate}
\item $\displaystyle \g=\cart + \sum_{\alpha\in\Phi}\g_\alpha$ (direct sum).
\item $\dim \g_\alpha = 1 \mbox{ for each } \alpha \in \Phi$.
\item Let $\alpha, \ \beta$ be two roots such that $\alpha+\beta\neq 0$.  Then $\g_\alpha$ and $\g_\beta$ are orthogonal under $B$.
\item The restriction of $B$ to $\cart \times \cart$ is nondegenerate.  For each $\alpha\in\carts$ there exists a unique element $H_\alpha\in\cart$ such that $$B(H,H_\alpha)=\alpha(H) \ \ \mbox{for all}\  H\in\cart.$$
\item  If $\alpha\in\Phi$, then $-\alpha\in\Phi$ and $$[\g_\alpha,\g_{-\alpha}]=\mathbb{C}H_\alpha, \ \ \alpha(H_\alpha)\neq 0.$$
\end{enumerate}
\label{roots}
\end{theorem}

We call the elements $H_\alpha\in\cart$ root vectors.  Let $\displaystyle \tau_\alpha=\frac{2H_\alpha}{B(H_\alpha,H_\alpha)}$ for $\alpha\in\Phi$,  then 
$\beta(\tau_\alpha)\in\mathbb{Z}$ for all $\alpha, \beta\in\Phi$.  We put $(\lambda, \mu) = B(H_\lambda,H_\mu)$ for $\lambda, \ \mu\in \carts$.

Let $\alpha\in\Phi$ and $\beta$ be any root.  Then by \cite[Theorem 4.3]{hel2} $[\g_\alpha,\g_\beta]=\g_{\alpha+\beta}$.

\begin{theorem}\cite[Theorem 4.4]{hel2}   Let $\cartr=\sum_{\alpha\in\Phi}\mathbb{R}H_\alpha$.  Then
\begin{enumerate}
\item $B$ is real and strictly positive definite on $\cartr\times\cartr$.
\item $\cart=\cartr\oplus i\cartr$.
\end{enumerate}
\end{theorem}

Recall that any element $\alpha\in\Phi$ determines a reflection $\sigma_\alpha\in\mbox{GL}(\carts)$ given by  $\displaystyle \sigma_\alpha(\beta)=\beta-\frac{2(\beta,\alpha)}{(\alpha,\alpha)}\alpha$. Then the Weyl group $\weyl$ of $\Phi$ is defined to be the subgroup of $GL(\carts)$ generated by the reflections $\sigma_\alpha$ for $\alpha\in\Phi$.  The following result gives that each root is conjugate to a simple root under action of the Weyl group.

\begin{theorem}\cite[Theorem 10.3.c]{humph}  Let $\Delta$ be a base of $\Phi$.  If $\alpha$ is any root, there 
exists $\sigma\in\mathcal{W}$ such that $\sigma(\alpha)\in\Delta$.\end{theorem}

We reduce to irreducible roots later, so we introduce the concept here.
If the set of nonzero roots $\Phi$ cannot be partitioned into the union of 2 proper orthogonal subsets, 
then $\Phi$ is called irreducible.  For any irreducible set of roots, at most 2 root lengths occur with roots of the same length conjugate under $\weyl$.
If $\Phi$ has two distinct root lengths they are described as long and short roots.  
In the case of only one root length, all roots are said to be long.  Any irreducible root system corresponds to one of the following types of classical Lie algebras: $A_n (n\geq 1), \ B_n (n\geq 2), \ C_n (n\geq 3), \ D_n
(n\geq 4), \ E_6, \ E_7,\ E_8, \ F_4,\ G_2$.  These results are proven in \cite{humph}.

Lie algebras of type $A_n, \ D_n, \ E_6, \ E_7, \mbox{ and } E_8$ have only long roots.  Those of type $B_n, \ C_n, \ F_4, \mbox{ and } G_2$ have both long and short roots.

The following results of Humphreys allow us to reduce our consideration of semisimple Lie algebras to simple Lie algebras.  Further discussion and proofs are found in Section 14.1 of \cite{humph}.

\begin{proposition}\cite[Proposition 14.1]{humph} Let $\g$ be a simple Lie algebra, $\cart$ a Cartan subalgebra and $\Phi$ the set of roots of $\g$ relative to $\cart$.  Then $\Phi$ is an irreducible root system.\label{humph14.1}
\end{proposition}

\begin{corollary}\cite[Corollary 14.1]{humph}  Let $\g$ be a semisimple Lie algebra with Cartan subalgebra $\cart$ and root system $\Phi$.  If  $\g=\g_1\oplus\cdots\oplus \g_t$ is the decomposition of $\g$ into simple ideals, then $\cart_i=\cart\cap \g_i$ is a Cartan subalgebra of $\g_i$, and the corresponding (irreducible) root system $\Phi_i$ may be regarded canonically as a subsystem of $\Phi$ in such a way that $\Phi=\Phi_1\cup\cdots\cup \Phi_t$ is the decomposition of $\Phi$ into its irreducible components.\label{cor14.1}
\end{corollary}

\subsection{Chevalley Basis}
\label{sec2.2}

We define a Chevalley basis for $\g$ a complex semisimple Lie algebra with root space decomposition
$\g=\cart+\sum_{\alpha\in\Phi}\g_\alpha$.  Fix a Cartan subalgebra $\cart$ of $\g$.  As above, let $H_\alpha$, $\alpha\in\Phi$ be the root vectors.   Let $\{X_\alpha\in\g_\alpha\}$ be elements that 
satisfy the following.
\begin{enumerate}
\item $[X_\alpha,X_{-\alpha}]=\tau_\alpha$ for all $\alpha\in\Phi$.
\item If $\alpha,\beta, \alpha+\beta\in\Phi$, then $[X_\alpha,X_\beta]=c_{\alpha\beta}X_{\alpha+\beta}$
where $c_{\alpha\beta}=-c_{-\alpha,-\beta}$.
\end{enumerate}
\begin{definition}
The set  $\chev=\{\tau_\alpha| \alpha \in \Delta; \ X_\beta| \beta\in\Phi\}$ is a \emph{Chevalley basis}
of $\g$, determined by $\cart$.\label{chevellay}
\end{definition}
Chevalley bases exist; see for example Section 25.2 of \cite{humph}.  Also, it follows from the definition of a Chevalley basis that the structure constants $\{C^\alpha_{ij}\}$  all lie in $\mathbb{Z}$, where $\{C^\alpha_{ij}\}$ are defined by $\displaystyle [X_i,X_j]=\sum_{\alpha\in\Phi} C^\alpha_{ij} X_\alpha$.

\subsection{Compact Chevalley Basis}
\label{sec4.3}

\begin{definition}For a Chevalley basis 
$\chev=\{\tau_\alpha|\alpha\in\Delta; \ X_\beta| \beta\in\Phi\}$ of $\g$  define
$\mathcal{C}_0=\{i\tau_\alpha|\alpha\in\Delta; A_\beta,B_\beta| \beta\in\Phi\}$ to be a \emph{compact Chevalley basis of $\goo$}  where 
 $A_\beta=X_\beta-X_{-\beta}$ and $B_\beta=i(X_\beta+X_{-\beta})$.   \end{definition}
It is known that $\mathbb{R}-\mbox{span}(\mathcal{C}_0)$ is a compact real form for $\g$, and one may choose the Cartan subalgebra $\cart$ so that $\goo=\mathbb{R}-\mbox{span}(\mathcal{C}_0)$ (see Theorem 6.3 of Ch. III of \cite{hel2}).  
If $\carto=\mathbb{R}-\mbox{span}\{i\tau_\alpha|\alpha\in\Delta\}$, then $\carto$ is a maximal abelian subalgebra of $\goo$, and  $\cartoc=\cart$.

Since the structure constants for a Chevalley basis are integers, it follows that the structure constants for a compact Chevalley basis $\mathcal{C}_0$ of $\goo$ also are integers.

\subsection{Chevalley rational structure}
\label{sec4.4}

Let $B(\mathbb{Q},\chev_0)$ be the set of bases $\mathcal{B}$ of $U$ such that $d\rho(\chev_0)$ leaves invariant $\mathbb{Q}-\mbox{span}(\mathcal{B})$.  By \cite{ragart}, the set $B(\mathbb{Q},\mathcal{C}_0)$ is nonempty.  For $\mathcal{B}\in B(\mathbb{Q},\mathcal{C}_0)$ let $\n_{\mathbb{Q}}=\mathbb{Q}-\mbox{span}\{\mathcal{B}\cup \mathcal{C}_0\}$.  One can show that for a suitably chosen $\goo-$invariant inner product on $U$, the structure constants of $\n=U\oplus\goo$ are rational (cf. \cite{ebpp}, Section 9 Proposition).  We call $\n_{\mathbb{Q}}$ a Chevalley rational structure for $\n$.

\begin{proposition}  Let $Z\in\mathbb{Q}-\mbox{span}(\chevo)$.  Then $\ker d\rho(Z)$ is a rational subspace of $U$ with respect to the Chevalley rational structure $\n_\mathbb{Q}$. \label{kerz}\end{proposition}

\begin{proof}
By the discussion above there exists  a basis $\mathcal{B}$ of $U$ such that $d\rho(Z)$ leaves invariant $\mathbb{Q}-\mbox{span}(\mathcal{B})$.  We identify $d\rho(Z)$ with its rational $q\times q$ matrix $Z^*$ defined by  the basis $\mathcal{B}$, where $q=\dim U$.

To find $\ker d\rho(Z)$, we consider the nullspace of the rational matrix $Z^*$.  To prove that $\ker d\rho(Z)$ has a basis of rational vectors it suffices to show that the dimension of the nullspace over $\mathbb{Q}$ is the same as the dimension of the nullspace over $\mathbb{R}$. 

Recall that $\mbox{rank}_\mathbb{Q}Z^*$ is the largest integer $r$ such that some $r\times r$ minor of $Z^*$ has nonzero determinant.  However this is the same as $\mbox{rank}_\mathbb{R}Z^*$ since all minor determinants of $Z^*$ are rational numbers.  Therefore $\dim_\mathbb{Q}\ker Z^*=q-\mbox{rank}_\mathbb{Q}Z^*=q-\mbox{rank}_\mathbb{R}Z^*=\dim_\mathbb{R}\ker Z^*$. \qed
\end{proof}

\subsection{The Basics of Weights}
\label{sec3.1}

We continue to let $\Phi$ denote the set of roots of a semisimple Lie algebra $\g$.   
Let $\Lambda=\{\lambda\in\mbox{Hom}(\cart,\mathbb{C}): \langle \lambda,\alpha\rangle \in\mathbb{Z} \ \forall \ \alpha\in\Phi\}$ be the set of abstract weights associated to $\cart$ a Cartan subalgebra of $\g$.  

Let $V$ be a finite dimensional complex $\g-$module.
With respect to $\cart$, a fixed Cartan subalgebra of $\g$, we consider the sets $V_\lambda=\{v\in V| Hv=\lambda(H)v \ \forall \ H \in\cart\}$.  If $V_\lambda \neq \{0\}$, then $V_\lambda$ is a weight space and $\lambda$ a weight of $\cart$ on $V$ or simply a weight of $V$.  We denote this set of weights determined by $V$ as $\Lambda(V)\subseteq\cart^*=\mbox{Hom}(\cart,\mathbb{C})$.

As in \cite{humph} we define  $\langle \alpha, \beta\rangle=\frac{2(\alpha,\beta)}{(\beta,\beta)}$ where $(\alpha,\beta)$ is any symmetric bilinear form on $\cart^*$ invariant under the Weyl group $\weyl$.  Note that $\langle \alpha,\beta\rangle$ is linear in the first variable, but not in the second.  It is known that $\langle \lambda, \alpha \rangle\in\mathbb{Z}$ for all $\linl$ and $\alpha\in\Phi$.  
For a fixed base $\Delta=\{\aone,\dots,\an\}$ of $\cart$, a weight $\linl$ is called dominant if $0\leq \langle \lambda,\alpha_i\rangle\in\mathbb{Z}$ for $1\leq i\leq n$.  We denote the set of all dominant weights by $\Lambda^+$.
We define fundamental dominant weights $\{\omega_1,\dots,\omega_n\}$ dual to the base $\Delta$ by requiring $\langle \omega_i,\alpha_j\rangle=\delta_{ij}$.

\begin{lemma}\cite[Lemma 13.2.A]{humph}  Each weight is conjugate under the action of the Weyl group $\weyl$ to one and only one 
dominant weight.\label{domweight}
\end{lemma}

Let $\Pi$ be a subset of $\Lambda$.  We call the set $\Pi$ saturated if for all $\lambda\in\Pi$ and
all $\alpha\in\Phi$, $\lambda-m\alpha$ is also in $\Pi$ for all $m$ between zero and $\langle \lambda, \alpha \rangle$.
A partial ordering of the set $\Lambda$ exists; we say that $\lambda \succ \mu$ if  $\lambda - \mu$ is a sum of positive roots and $\lambda\succeq\mu$ if $\lambda=\mu$ or $\lambda\succ\mu$.  This leads to a natural definition of a highest weight of $\Pi$ to be that weight $\lambda\in\Pi$ such that 
$\lambda \curvegreat \mu$ for all other $\mu\in\Pi$.

\begin{lemma} \cite[Lemma 13.4.A]{humph}  If $\lambda\in\Lambda^+$, then a saturated set of weights having highest weight $\lambda$ 
must be finite.\label{humphA}\end{lemma}

\begin{lemma}\cite[Lemma 13.4.B]{humph}   Let $\Pi$ be saturated with highest weight $\lambda\in\Lambda^+$.  If
$\mu\in\Lambda^+$ and $\lambda\succeq \mu$, then $\mu\in\Pi$.\label{humphB}\end{lemma}

From Lemmas \ref{domweight} and \ref{humphB} we obtain a complete description of a saturated set $\Pi\subseteq\Lambda$ of highest weight $\lambda$.  The set $\Pi$ consists of all
dominant weights $\mu$ such that $\mu\preceq\lambda$ and the conjugates
of $\mu$ under the Weyl group $\weyl$.

\subsection{Weight Spaces}
\label{sec3.2}

Fix a Cartan subalgebra $\cart$ of $\g$.  Let $V$ be a complex $\g-$module and let $\Lambda(V)\subseteq \mbox{Hom}(\cart,\mathbb{C})$ be the set of weights determined by $\cart$ and $V$.  Then $\displaystyle V=V_0\oplus \sum_{\lambda\in\Lambda }V_\lambda$ is the weight space 
decomposition of $V$.   

\begin{lemma}\cite[Lemma 20.1]{humph}   Let $V$ be an arbitrary $\g-$module.  Then $\g_\alpha$ maps $V_\lambda$ into
$V_{\lambda+\alpha}$.
\label{humph107}
\end{lemma}

\begin{theorem} Let $\lambda\in\Lambda^+\subseteq\cart^*.$  Then there exists an irreducible finite dimensional $\g-$module $V(\lambda)$, unique up to isomorphism, with highest weight $\lambda$.  The map $\lambda\rightarrow V(\lambda)$ is a bijection between $\Lambda^+$ and the isomorphism classes of finite dimensional irreducible $\g-$modules.
\end{theorem}

A proof can be found in \cite{humph}, p 113.

\subsection{Closed Geodesics and the Set $L'$}
\label{sec3.3}

If $\goo$ is a compact semisimple Lie algebra and $U$ is a finite dimensional real $\goo-$module, then $\n=U\oplus\goo$ admits a metric, 2-step nilpotent Lie algebra structure defined in \ref{sec4.1}.  Let $N$ denote the simply connected 2-step nilpotent Lie group with Lie algebra $\n$ and left invariant metric $\ip$ arising from the inner product $\ip$ on $\n$.  The discussions in \ref{sec1.3} and \ref{sec4.3} show that $N$ admits a lattice $\Gamma$.  To study the density of closed geodesics  (DCG) on $\ngam$ we first consider the complex $\g=\goc-$module $V=U^\mathbb{C}$, and we assume that the zero weight space $V_0$ is nontrivial for a Cartan subalgebra $\cart$ of $\g$.  If $V_0=\{0\}$, then the DCG property is completely settled by the main result of \cite{leepark}.  For convenience we assume that $\cart=\cartoc$, where $\carto$ is a maximal abelian subspace of $\goo$.

Our method in studying the closed geodesics of $\ngam$ requires that every root of $\cart$ be a weight of $V=U^\mathbb{C}$ of complex multiplicity at least two.  This is explained in Section \ref{chap6}.  Furthermore we may reduce to the case $\goo$ simple (cf. Section \ref{chap6}).  In this case $\g=\goc$ is a complex simple Lie algebra.

Let $L'$ be the set of dominant weights such that the associated irreducible $\g-$module $V(\lambda)$ fails to satisfy this multiplicity condition on the roots of $\cart$.  Using work of \cite{bz} we determine the list $L'$ of Table \ref{lprime}, first finding those highest weights  for which not all roots are weights and then considering the multiplicity condition.  In Section \ref{chap6} we explain the relevance of the list $L'$ to the solution of the DCG problem for $\ngam$.

\begin{table}[h]
\caption{The list $L'$}
\centering
\label{lprime}
\begin{tabular}{|c|c|} \hline
Lie algebra type & dominant highest weight \\ \hline
$A_n$ & $\alpha_1+\alpha_2+\cdots+\alpha_n$\\
 	& $k(n\alpha_1+(n-1)\alpha_2+\cdots +\alpha_n)$, $k\in\mathbb{N}$, $k> 1$\\ 
	& $\alpha_1+2\alpha_2$, $2\alpha_1+\alpha_2$ $(n=2)$  \\
	& $\alpha_1+2\alpha_2+\alpha_3$ $(n=3)$  \\ \hline
$B_n$ $(n\geq 2)$ & $\alpha_1+\alpha_2+\cdots +\alpha_n$ \\ 
 	& $\alpha_1+2\alpha_2\cdots +2\alpha_n$ \\
	& $\alpha_1+2\alpha_2+m_3\alpha_3$, $m_3\geq 3$ $(n=3)$\\
	& $2\alpha_1+2\alpha_2\cdots +2\alpha_n$ $(n\geq2)$\\ 
	& $\alpha_1+2\alpha_2+3\alpha_3+3\alpha_4+\cdots +3\alpha_n$ $(n\geq3)$\\ \hline
$C_n$ $(n\geq 3)$ & $\alpha_1+2\alpha_2+2\alpha_3+\cdots+2\alpha_{n-1}+\alpha_n$\\
	& $2\alpha_1+2\alpha_2+2\alpha_3+\cdots+2\alpha_{n-1}+\alpha_n$\\
	& $\alpha_1+2\alpha_2+\cdots+(n-1)\alpha_{n-1}+[\frac{n}{2}]\alpha_n$ \\ 
	& $\alpha_1+2\alpha_2+\cdots+2N\alpha_{2N}+\cdots+2N\alpha_{n-1}+N\alpha_n$ \\ \hline
$D_n$ $(n\geq 4)$  & $\alpha_1+2\alpha_2+\cdots +2\alpha_{n-2}+\alpha_{n-1}+\alpha_n$\\
	&$2\alpha_1+2\alpha_2+\cdots+2\alpha_{n-2}+\alpha_{n-1}+\alpha_n$\\ 
&  $\alpha_1+2\alpha_2+2\alpha_3+\alpha_4$, $\alpha_1+2\alpha_2+\alpha_3+2\alpha_4$ $(n=4)$\\ \hline
$E_6$	& $\alpha_1+2\alpha_2+2\alpha_3+3\alpha_4+2\alpha_5+\alpha_6$\\ \hline
$E_7$	& $2\alpha_1+2\alpha_2+3\alpha_3+4\alpha_4+3\alpha_5+2\alpha_6+\alpha_7$\\ \hline
$E_8$	& $2\alpha_1+3\alpha_2+4\alpha_3+6\alpha_4+5\alpha_5+4\alpha_6+3\alpha_7+2\alpha_8$\\ \hline
$F_4$ & $2\alpha_1+3\alpha_2+4\alpha_3+2\alpha_4$ \\ 
	& $\alpha_1+2\alpha_2+3\alpha_3+2\alpha_4$ \\ \hline
$G_2$ & $2\alpha_1+\alpha_2$ \\ 
	& $3\alpha_1+2\alpha_2$ \\
	& $4\alpha_1+2\alpha_2$ \\ \hline
\end{tabular}
\end{table}

\subsection{Density of Resonant and Rational Vectors}
\label{sec4.5}

Recall from \ref{sec1.3} that a rational basis for $\n=U\oplus\goo$ is one for which all the structure constants are rational.  We say a vector is rational if it lies in the rational span of a rational basis.  In this section we now fix a rational basis $\mathcal{B}\cup \mathcal{C}_0$, where $\mathcal{B}\in B(\mathbb{Q},\mathcal{C}_0)$.

\begin{proposition} Let $\n=U\oplus\g_0$.  The set of rational, resonant vectors in $\goo$ is dense
in $\goo$.\label{ratlres}\end{proposition}

To prove this proposition, we need the following result.  We are interested in (2) of the  theorem, which is a consequence of (1) and a result of Sansuc \cite[Corollary 3.5 (iii)]{sansuc}.  For $\mathcal{B}\in B(\mathbb{Q},\chev_0)$,
define the subgroup $G_{\mathcal{B},\mathbb{Q}}=\{g\in G_0| \rho(g) \mbox{ leaves invariant }\mathbb{Q}-\mbox{span}\{\mathcal{B}\}\}$.
\begin{theorem}\cite[Section 5, Theorem A]{eb2}  Let $\rho:G_0\rightarrow GL(U)$ be as above, and for a basis $\mathcal{B}$ of $U$, let $\rho_\mathcal{B}: G_0 \rightarrow GL(n,\mathbb{R})$ denote the corresponding Lie group homomorphism.
\begin{enumerate}
\item If $\mathcal{B}\in B(\mathbb{Q},\chevo)$, then $\rho_\mathcal{B}(G)$ is an affine algebraic group defined over $\mathbb{Q}$.
\item If $\mathcal{B}\in\ B(\mathbb{Q},\chevo)$, then $G_{\mathcal{B},\mathbb{Q}}$  is dense in $G$ in the Lie topology.
\end{enumerate}
\label{ebthmA}
\end{theorem}
\begin{proof}[Proof of Proposition \ref{ratlres}]
Let $G_0$ be the compact, simply connected Lie group with Lie algebra $\goo$ and let $\rho: G_0 \rightarrow GL(U)$ be the  real representation corresponding to $d\rho: \goo\rightarrow \mbox{End}(U)$.  Let $d\tilde{\rho}: \g\rightarrow \mbox{End}(V)$ be the associated complex representation, where $\g=\goc$ and $V=U^{\mathbb{C}}$.

For the finite subset $\Lambda_\rho\subseteq \mbox{Hom}(\cart,\mathbb{C})$ of weights, we write the weight space decomposition $\displaystyle V=\sum_{\beta\in{\Lambda_\rho}}V_\beta$ where $d\tilde{\rho}(H)=\beta(H)Id$ on the subspace $V_\beta$ for each $H\in\cart$ and $\beta\in\Lambda_\rho$.  By standard representation theory, $d\tilde{\rho}(\tau_\alpha)$ has eigenvalues $\lambda(\tau_\alpha)\in\mathbb{Z}$ for all $\alpha\in\Phi$ and $\lambda\in\Lambda_\rho$.

For convenience we let $\tilde{\tau_\alpha}=i\tau_\alpha\in\mathcal{C}_0$.  The elements of the set $\{\tilde{\tau}_\alpha| \alpha\in\Delta\}$ commute, and thus we can find a common basis of eigenvectors in $U$ for the space $\cart_\mathbb{Q}=\mathbb{Q}-\mbox{span}\{\tilde{\tau}_\alpha| \alpha\in\Delta\}$.  Since $d\rho(i\tilde{\tau}_\alpha)=d\tilde{\rho}(-\tau_\alpha)$ and $d\tilde{\rho}(\tau_\alpha)$ has eigenvalues in $\mathbb{Z}$, $d\rho(\tilde{\tau}_\alpha)$ has eigenvalues in $i\mathbb{Z}$.  It follows that the elements of $\cart_{\mathbb{Q}}$ have eigenvalues in $i\mathbb{Q}$ and hence are resonant.  Since $\{\tilde{\tau}_\alpha| \alpha\in\Delta\}$  is a basis for $\carto$, $\mathbb{Q}-\mbox{span}\{\tilde{\tau}_\alpha| \alpha\in\Delta\}=\cart_\mathbb{Q}$ is dense in $\carto$.  Thus $\carto$ contains a dense set of rational, resonant vectors $\cart_\mathbb{Q}$.  Now we show that $\goo$ contains such a dense set, thus proving the proposition.

Let $\g_{\mathbb{Q}}=\mbox{Ad}(G_{\mathcal{B},\mathbb{Q}})(\cart_{\mathbb{Q}})$. We have shown that $\cart_{\mathbb{Q}}$ is dense in $\carto$, $G_{\mathcal{B},\mathbb{Q}}$ is dense in $G_0$, and $\mbox{Ad}(G_0)(\carto)=\goo$, thus it follows that $\g_{\mathbb{Q}}$ is dense in $\goo$.  We show that the elements of $d\rho(\cart_{\mathbb{Q}})$ have eigenvalues in $i\mathbb{Q}$ and hence are resonant.  If $X\in\g_{\mathbb{Q}}$, then $X=\mbox{Ad}(g)A$ for $g\in G_{\mathcal{B},\mathbb{Q}}$ and $A\in\cart_{\mathbb{Q}}$.  It follows that $d\rho(X)=\rho(g)d\rho(A)\rho(g)^{-1}$ has the same eigenvalues in $i\mathbb{Q}$ as $d\rho(A)$. Since both $d\rho(\cart_\mathbb{Q})$ and $\rho(G_{\mathcal{B},\mathbb{Q}})$ leave $\mathbb{Q}-\mbox{span}(\mathcal{B})$ invariant, all elements of $d\rho(\g_\mathbb{Q})$ leave $\mathbb{Q}-\mbox{span}(\mathcal{B})$ invariant.  Thus $\g_\mathbb{Q}$ is a set of rational, resonant vectors dense in $\goo$.\qed \end{proof}

\section{First Hit Map}
\label{chap5}

Let $\{\n=U\oplus \goo, \ \ip\}$ be a metric 2-step nilpotent Lie algebra with associated simply connected Lie group $N$, where
 $\goo$ is an admissible semisimple compact real Lie algebra as defined below and 
$U$ is a real $\goo$-module.  We shall assume further that $\{0\}=\ker(\goo)=\{u\in U|Zu=0 \ \forall Z\in\goo\}$.  Under this condition $\goo=\z$, the center of $\n=U\oplus\goo$.

Note that $U$ is a $G_0-$module, where $G_0$ is the compact, simply connected Lie group with Lie algebra $\goo$.  Let $\carto$ be a maximal abelian subalgebra of $\goo$
and $\cart=\cartoc$ be a Cartan subalgebra of $\g=\goc$.   Let $\chevo$ be a compact Chevalley basis of $\goo$.  Let $\mathcal{B}$ be a basis of $U$ such that $\chevo$ leaves $\mathbb{Q}-\mbox{span}(\mathcal{B})$ invariant and let $\n_{\mathbb{Q}}=\mathbb{Q}-\mbox{span}\{\mathcal{B}\cup\chevo\}$ be the Chevalley rational structure for $\n$.

Let $\Lambda\subseteq \mbox{Hom}(\cart,\mathbb{C})$ be the set of weights associated to the real $\goo-$module $U$ as described in Section \ref{chap4}.  From Proposition \ref{eb4.2}, if $\linl$, then $-\linl$.  Moreover $\lambda(Z)\in i\mathbb{R}$ for all $Z\in\carto$, $\linl$ since the elements of $\goo$ are skew symmetric on $U$.

\subsection{Definition of the First Hit Map}
\label{sec5.1}

\begin{definition}A vector $Z\in\carto$ is \emph{super regular} if $\lambda(Z)\neq 0$ for all nonzero $\linl$ and $\lambda(Z)^2=\mu(Z)^2$ for $\lambda,\mu\in\Lambda$ implies that $\mu=\pm\lambda$.  For $Z\in\goo$ we say that $Z$ is \emph{super regular} if there exists $g\in G_0$ such that $\mbox{Ad}(g)Z\in\carto$ and $\mbox{Ad}(g)Z$ is super regular.  
\end{definition}
One can show that the element $\mbox{Ad}(g)$ is uniquely determined and the definition of a super regular vector $Z\in\goo$ is independent of the choice of maximal abelian subalgebra $\carto$.  
The set of super regular vectors $Z$ in $\carto$ is a dense open subset of $\carto$ since it is the complement of a finite union of hyperplanes in $\carto$.  It follows from $\goo=\mbox{Ad}(G_0)(\carto)$ that the super regular vectors in $\goo$ form a dense open subset of $\goo$.  In particular, the super regular vectors in $\cart_{\mathbb{Q}}$ and $\g_{\mathbb{Q}}=\mbox{Ad}(G_{\mathcal{B},\mathbb{Q}})(\cart_{\mathbb{Q}})$ are dense in $\cart_{\mathbb{Q}}$ and $\g_{\mathbb{Q}}$ respectively.  For $Z$ super regular, note that the eigenvalues $\{0,\lambda(Z)| \linl\}$ are all distinct.  
In the following, we assume that $Z$ is super regular. We let $j:\goo\rightarrow \mbox{End}(U)$ denote the representation.

 Let $\gamma(t)$
be a geodesic in $N$ with $\gamma(0)=e$ and $\gamma'(0)=X+ \alpha Z$, where $X\in U$, $Z\in\goo$ and $0\neq \alpha\in\mathbb{R}$.
We consider first the case that $Z\in\carto$ and later reduce the general case to this one.  Write $U=U_0\oplus U_1$ where $U_0=\ker j(Z)$ and $U_1$ is the orthogonal complement.
Then $\displaystyle U_1=\bigoplus_{\linl}U_\lambda$, where $U_\lambda$ is the eigenspace of $j(Z)^2$ corresponding to $-a_\lambda^2=-\lambda(Z)^2$.  For each $X\in U$, $X=X_0+X_1$, where
$X_i\in U_i$ for $i=0,1$ and $\displaystyle X_1=\sum_{\linl}X_\lambda$, $X_\lambda\in U_\lambda$.
Thus we have $\displaystyle X=X_0+\sum_{\lambda\in\Lambda}X_\lambda$.

Recall that $\z=\goo$ is the center of $\n=U\oplus\goo$ since $\ker(\goo)=\{0\}$.  Let $\z_{res}$ denote the set of resonant vectors in $\z$.  For $Z\in\z_{res}$, we define the sets $\n_Z$ and $\w_Z$ as follows:
\begin{eqnarray*}
\n_Z&=&\{X+\alpha Z\in\n| X\in U,\ X \mbox{ has a nonzero component in }\ker j(Z) \\ & & \mbox{ and } 0\neq \alpha\in\mathbb{R}\}\\
\w_Z&=&\mbox{span}\{\z,\ker j(Z)\}=\mbox{span}\{\goo,\ker j(Z)\}
\end{eqnarray*}

\begin{remark} \label{remark4.1} If $Z\in\mathbb{Q}-\mbox{span}(\chevo)$, then $\w_Z$ is a rational subalgebra with respect to any Chevalley rational structure on $\n$ by Proposition \ref{kerz} and the fact that $\goo$ is a rational subspace of $\n$ with rational basis $\chevo$.
\end{remark}

 For any $\xi=X+\alpha Z\in\n_Z$, let 
$\gamma_\xi(t)$ be the geodesic in N with $\gamma_\xi(0)=e$ and 
$\gamma_\xi'(0)=\xi$. We define the first hit map $F_Z: \n_Z\rightarrow \w_Z$ as
$$F_Z(\xi)=\log (\gamma_\xi(\omega)) \mbox{ where }\omega>0 \mbox{ is the first value for which } e^{\omega j(\alpha Z)}=Id.$$

Write $\gamma_\xi(t)=\exp(X(t)+Z(t))$, where $X(t)\in U$ and $Z(t)\in\goo$.   Eberlein \cite{eb1} has shown that the geodesic equations for $\gamma_\xi(t)$ are given by
\begin{eqnarray*}
X(t)&=&tX_0+(e^{tj(\alpha Z)}-Id)(j(\alpha Z)^{-1}X_1)\\
Z(t)&=&t\tilde{Z_1}(t)+\tilde{Z_2}(t)\end{eqnarray*}
where
\begin{eqnarray*}
\tilde{Z_1}(t)&=&\alpha Z +\frac{1}{2}[X_0,(e^{tj(\alpha Z)}+Id)j(\alpha Z)^{-1}X_1]+\frac{1}{2}\sum_{\lambda\in\Lambda}
[j(\alpha Z)^{-1}X_\lambda,X_\lambda]\\
\tilde{Z_2}(t)&=&[X_0,(Id-\ej)j(\alpha Z)^{-2}X_1]+\frac{1}{2}[\ej\jinv X_1,\jinv X_1]\\
& & -\frac{1}{2}\sum_{\lambda\neq \mu \in \Lambda}\left( \frac{1}{a_\mu^2-a_\lambda^2}\right)
([\ej j(\alpha Z)X_\lambda, \ej \jinv X_\mu]-[\ej X_\lambda, \ej X_\mu])\\
& & +\frac{1}{2}\sum_{\lambda\neq \mu \in \Lambda}\left( \frac{1}{a_\mu^2-a_\lambda^2}\right)
([j(\alpha Z)X_\lambda, \jinv X_\mu]-[X_\lambda,X_\mu])\\
\end{eqnarray*}

Here $\jinv$ denotes the inverse of $j(\alpha Z)$ on $\displaystyle U_0^\perp=U_1=\bigoplus_{\linl}U_\lambda$.  Thus $\jinv=\frac{1}{\alpha}j(Z)^{-1}$.
It follows from the geodesic equations that $F_Z(\xi)\in\w_Z$ since $X(\omega)=\omega X_0\in U_0$ and $Z(\omega)\in\goo$.

\subsection{Properties of the First Hit Map}
\label{sec5.2}

The proofs of the next two results parallel Lemmas 14 and 15 of \cite{demeyer} and the proofs.  We omit the details.

\begin{lemma} Let $\alpha$ be any nonzero real number and $Z$ any nonzero element of $\z_{res}$.  Then, 
\begin{enumerate}
\item $\n_Z=\n_{\alpha Z}$ and $\w_Z=\w_{\alpha Z}$,
\item $F_Z=F_{\alpha Z}$, and
\item $F_Z(\alpha\xi)=F_Z(\xi)$ for all $\xi\in\n_Z$.
\end{enumerate}
\label{demeyer14}
\end{lemma}

Recall from Theorem \ref{laur} that for each $g\in G_0$, there exists a map $I_g\in\mbox{Aut}(N)\cap I(N)$ such that $dI_g$ acts on  $\n=U\oplus\goo$  as an automorphism and an isometry by $(\rho(g),\mbox{Ad}(g))$.  The induced map  $dI_g: \n \rightarrow \n$ preserves the center $\z$ since it is an automorphism and preserves the orthogonal complement $\z^{\perp}$ since it is an isometry.%

\begin{lemma} Let $g$ be any element of $G_0$.  Then $F_{Ad(g)Z}=dI_g\circ F_Z\circ(dI_g)^{-1}$ for all nonzero $Z\in\z_{res}$.
\label{demeyer15}
\end{lemma}

These two results allow us to reduce the question of maximal rank of the first hit map for $Z\in\goo$ to $Z\in\carto$ by the following corollary and the fact that $\mbox{Ad}(G_0)(\carto)=\goo$.

\begin{corollary}  Let $Z\in\goo$ be resonant.  Then $F_Z:\n_Z\rightarrow \w_Z$ has maximal rank if and only if $F_{Ad(g)Z}:\n_{Ad(g)Z}\rightarrow \w_{Ad(g)Z}$ has maximal rank for any $g\in G_0$.\label{nicez}
\end{corollary}

\subsection{The First Hit Map for $\n=U\oplus\goo$}
\label{sec5.3}

Let $Z$ be a nonzero element of $\z_{res}$.  For $\n=U\oplus \goo$ and $\xi=X+\alpha Z\in\n_Z$, we show 

\begin{equation}
F_Z(\xi)=\omega\{X_0+\alpha Z+\frac{1}{\alpha}[X_0,j(Z)^{-1} X_1]+\frac{1}{2\alpha}\sum_{\linl}[j(Z)^{-1} X_\lambda, X_\lambda]\}.
\label{fsthit}
\end{equation}

We calculate $X(\omega)$ and $Z(\omega)$ recalling 
that $\omega$ is chosen so that $e^{\omega j(\alpha Z)}=Id$.  Thus the terms are greatly simplified as follows.
\begin{eqnarray*}
X(\omega)&=&\omega X_0+(e^{\omega j(\alpha Z)}-Id)(j(\alpha Z)^{-1}X_1)\\
&=& \omega X_0\\
Z(\omega)&=&\omega \tilde{Z_1}(\omega)+\tilde{Z_2}(\omega) 
\end{eqnarray*}
where
\begin{eqnarray*}
\tilde{Z_1}(\omega)&=&\alpha Z +\frac{1}{2}[X_0,(e^{\omega j(\alpha Z)}+Id)j(\alpha Z)^{-1}X_1]+\frac{1}{2}\sum_{\lambda\in\Lambda}
[j(\alpha Z)^{-1}X_\lambda,X_\lambda]\\
&=&\alpha Z +[X_0,j(\alpha Z)^{-1}X_1]+\frac{1}{2}\sum_{\lambda\in\Lambda}
[j(\alpha Z)^{-1}X_\lambda,X_\lambda]\\
&=&\alpha Z +\frac{1}{\alpha}[X_0,j(Z)^{-1}X_1]+\frac{1}{2\alpha}\sum_{\lambda\in\Lambda}
[j(Z)^{-1}X_\lambda,X_\lambda]\\
\tilde{Z_2}(\omega)&=&[X_0,(Id-\ejo)j( \alpha Z)^{-2}X_1]+\frac{1}{2}[\ejo\jinv X_1,\jinv X_1]\\
& & -\frac{1}{2}\sum_{\lambda\neq \mu \in \Lambda}\left( \frac{1}{a_\mu^2-a_\lambda^2}\right)
([\ejo j(\alpha Z)X_\lambda, \ejo \jinv X_\mu]-[\ejo X_\lambda, \ejo X_\mu])\\
& & +\frac{1}{2}\sum_{\lambda\neq \mu \in \Lambda}\left(\frac{1}{a_\mu^2-a_\lambda^2}\right)
([j(\alpha Z)X_\lambda, \jinv X_\mu]-[X_\lambda,X_\mu])\\
&=&\frac{1}{2}[\jinv X_1,\jinv X_1]\\
& & -\frac{1}{2}\sum_{\lambda\neq \mu \in \Lambda}\left( \frac{1}{a_\mu^2-a_\lambda^2}\right)
([j(\alpha Z)X_\lambda, \jinv X_\mu]-[ X_\lambda, X_\mu])\\
& & +\frac{1}{2}\sum_{\lambda\neq \mu \in \Lambda}\left( \frac{1}{a_\mu^2-a_\lambda^2}\right)
([j(\alpha Z)X_\lambda, \jinv X_\mu]-[X_\lambda,X_\mu])\\
&=& 0 \mbox{ by properties of the bracket.}
\end{eqnarray*}
This gives
\begin{eqnarray*}
X(\omega)&=&\omega X_0\\
Z(\omega)&=&\omega\left(\alpha Z+\frac{1}{\alpha }[X_0,j(Z)^{-1}X_1]+\frac{1}{2\alpha }\sum_{\lambda\in\Lambda}[j(Z)^{-1}X_\lambda,X_\lambda]\right).
\end{eqnarray*}
Thus since $F_Z(\xi)=X(\omega)+Z(\omega)$ by the definition of the first hit map, equation (\ref{fsthit}) is verified.

\subsection{Maximal Rank of the First Hit Map}  
\label{sec5.4}

In the terminology of \ref{sec4.4}, we recall that $\cart_{\mathbb{Q}}=\mathbb{Q}-\mbox{span}\{\tilde{\tau}_\alpha: \alpha\in\Delta\}$ and $\g_{\mathbb{Q}}=\mbox{Ad}(G_{\mathcal{B},\mathbb{Q}})(\cart_{\mathbb{Q}})$.  In \ref{sec4.4} we showed that $\g_{\mathbb{Q}}$ consists of resonant, rational vectors and that $\cart_{\mathbb{Q}}$ and $\g_{\mathbb{Q}}$ are dense in $\carto$ and $\goo$ respectively.  Let $\g_{\mathbb{Q},S}$ denote the super regular vectors in $\g_{\mathbb{Q}}$.  Recall that $\g_{\mathbb{Q},S}$ is dense in $\g_\mathbb{Q}$ and hence dense in $\goo$ by the discussion in \ref{sec5.1}.  We will show that $F_Z$ has maximal rank on a dense open subset of $\n_Z$ for every $Z\in\g_{\mathbb{Q},S}$.  By Corollary \ref{nicez}, it suffices to consider $Z\in\cart_{\mathbb{Q}}\subseteq\carto$ super regular.

For $Z\in\cart_{\mathbb{Q}}$ super regular, let $\xi=X+\alpha Z\in\n_Z$.  By the definition of $\n_Z$,  $0\neq \alpha\in\mathbb{R}$ and $\displaystyle X=X_0+\sum_{\linl}X_\lambda$, where $0\neq X_0\in U_0$ and $X_\lambda\in U_\lambda$ for every $\linl$.  Recall that $\Phi\subseteq \Lambda$ by hypothesis.  For $\beta\in\Phi$, we consider the maps $A_\beta,B_\beta:U_0\rightarrow U_\beta$, where $A_\beta=X_\beta-X_{-\beta}$ and $B_\beta=iX_\beta+iX_{-\beta}$ (c.f. \ref{sec4.3}).   Note that $\ker(A_\beta)=\ker(B_\beta)$ by Proposition 5.6 of \cite{ebpp}.  

Let \begin{eqnarray*}\n_Z^*&=&\{X+\alpha Z\in\n_Z| X_0\notin \ker(A_\beta) \mbox{ and } X_\beta\notin \mbox{span}\{A_\beta(X_0), B_\beta(X_0)\} \\ & & \mbox{ for all }\beta\in\Phi\}.\end{eqnarray*} It is evident that $\n_Z^*$ is a dense open subset of $\n_Z$.  

For $Z\in\g_{\mathbb{Q},S}$ we define $\n_Z^*$ as follows.  Since $Z$ is super regular, there exists a unique $g\in G_{\mathcal{B},\mathbb{Q}}$ such that $\mbox{Ad}(g)Z\in\cart_{\mathbb{Q}}$.  Now define $$\n_Z^*=\{\xi\in\n_Z|dI_g(\xi)\in \n_{Ad(g)Z}^*\},$$ where $\n_{Ad(g)Z}^*$ is defined as above.  In particular, $\n_Z^*$ is again a dense open subset of $\n_Z$.

\vspace{0.1in}

\begin{remark} \label{remark5.5.1} If $X+\alpha Z\in \n_Z^*$, then from Propositions 5.7 and 5.9 of \cite{ebpp}, we obtain the following results:
\begin{enumerate}
\item $\mbox{ad}(X_0):U_\beta\rightarrow \goo(\beta)$ is surjective
\item $\ker\mbox{ad} (X_0) \nsubseteq \ker \mbox{ad}(X_\beta)$ for every $\beta\in\Phi$
\end{enumerate}\end{remark}

\begin{proposition}  Let $Z\in\cart_{\mathbb{Q}}$ be super regular.  Then $F_Z:\n_Z\rightarrow \w_Z$ has maximal rank at every point of $\n_Z^*$.
\label{maxrank}
\end{proposition}

The next corollary follows directly from Corollary \ref{nicez}, Proposition  \ref{maxrank} and the definition of $\g_{\mathbb{Q}}$.

\begin{corollary} Let $Z\in\g_{\mathbb{Q}}$ be super regular.  Then $F_Z:\n_Z\rightarrow \w_Z$ has maximal rank at every point of $\n_Z^*$.
\label{maxrankg}
\end{corollary}

\begin{proof}[Proof of Proposition \ref{maxrank}]

Let $\xi=X+\alpha Z\in\n_Z^*$, where $Z\in\cart_{\mathbb{Q}}$ is super regular, and write $\displaystyle X=X_0+X_1=X_0+\sum_{\linl}X_\lambda\in U$.  Let $\xi(s)=\nu(s)+X+\alpha Z$,  where $\nu(s)=su_\lambda$ is of two possible types:
\begin{enumerate}
\item $\nu(s)=su_0, \ u_0\in\ker j(Z)$ or
\item $\nu(s)=su_\beta$, $u_\beta \in U_\beta, \ 0\neq \beta\in\Phi\subseteq \Lambda$.
\end{enumerate}
We will use the notation $F_Z(\xi(s))=F_Z(u_0(s))$ in the first case and $F_Z(\xi(s))=F_Z(u_\beta(s))$ in the second.  Note that $\displaystyle dF_Z(u_0)=\frac{d}{ds}\raisebox{-3mm}{\rule{0.1mm}{10mm}}_{\ s=0} F_Z(u_0(s))$ and $\displaystyle dF_Z(u_\beta)=\frac{d}{ds}\raisebox{-3mm}{\rule{0.1mm}{10mm}}_{\ s=0} F_Z(u_\beta(s))$.
\begin{proposition}The derivative of the first hit map at $\xi=X+\alpha Z$ in the above cases is
\begin{enumerate}
\item $\displaystyle dF_Z(u_0)=\omega\left(u_0+\frac{1}{\alpha }[u_0,\jjinv X_1]\right)$
\item $\displaystyle dF_Z(u_\beta)=\frac{\omega}{\alpha }\left([X_0,\jjinv u_\beta]+[\jjinv u_\beta,X_\beta]\right)$
\end{enumerate}\label{dFz}
\end{proposition}
\begin{proof}
In the first case we have 
\begin{eqnarray*}
F_Z(u_0(s))&=&\omega( su_0+X_0+ \alpha Z+\frac{1}{\alpha }[su_0+X_0,\jjinv X_1]\\
 & & +\frac{1}{2\alpha }\sum_{\linl}[\jjinv X_\lambda, X_\lambda])
\end{eqnarray*}

Differentiating with respect to $s$, 
\begin{eqnarray*}
dF_Z(u_0) &=&\frac{d}{ds}\raisebox{-3mm}{\rule{0.1mm}{10mm}}_{\ s=0} F_Z(u_0(s))\\
&=&\omega\left(u_0+\frac{1}{\alpha }[u_0,\jjinv X_1]\right).
\end{eqnarray*}

In the second case we have
\begin{eqnarray*}
F_Z(u_\beta(s))&=&\omega (X_0+\alpha Z+\frac{1}{\alpha }[X_0,\jjinv (su_\beta+X_1)]+\frac{1}{2\alpha }\sum_{\lambda\neq \beta}[\jjinv X_\lambda, X_\lambda]\\
& & +\frac{1}{2\alpha }[\jjinv (su_\beta +X_\beta),(su_\beta+X_\beta)]).
\end{eqnarray*}
Thus
\begin{eqnarray*}
\frac{d}{ds}F_Z(u_\beta(s))&=&\omega\left(\frac{1}{\alpha}[X_0,\jjinv u_\beta]+\frac{1}{2\alpha }\left([\jjinv u_\beta, su_\beta+X_\beta]+
[\jjinv (su_\beta+X_\beta),u_\beta]\right)\right)\\
dF_Z(u_\beta)&=&\frac{d}{ds}\raisebox{-3mm}{\rule{0.1mm}{10mm}}_{\ s=0} F_Z(u_\beta(s))  \\
&=&\frac{\omega}{\alpha} \left([X_0,\jjinv u_\beta]+ \frac{1}{2}([\jjinv u_\beta,X_\beta]+
[\jjinv X_\beta,u_\beta])\right)\\
&=&\frac{\omega}{\alpha}\left([X_0,\jjinv u_\beta]+ [\jjinv u_\beta,X_\beta] \right)\mbox{ by Lemma \ref{jbeta}.}
\end{eqnarray*}\qed
\end{proof}

\begin{lemma} For any $Z\in\goo$ and any $u_\beta,X_\beta\in U_\beta$  $$[\jjinv u_\beta,X_\beta]=[\jjinv X_\beta,u_\beta].$$ \label{jbeta}\end{lemma}

\begin{proof} If $\beta\in\Phi\subseteq\Lambda$ we know by Corollary  \ref{eb4.8} that
$[U_\beta,U_\beta]=\mathbb{R}-\mbox{span}\{\tilde{H}_\beta\}$, where $\tilde{H}_\beta$ denotes the weight vector in $\carto$ determined by $\beta$. 
Since $\jjinv X_\beta,$  $\jjinv u_\beta, \ X_\beta, \ u_\beta$ are all in $U_\beta$ by  \cite[Proposition 4.3]{ebpp}, both brackets
lie in $\mathbb{R}-\mbox{span}\{\tilde{H}_\beta\}=\langle \tilde{H}_\beta\rangle$.  It suffices to show that $[j(Z)^{-1}u_\beta,X_\beta]=[j(Z)^{-1}X_\beta,u_\beta]$ have the same inner product with $\tilde{H}_\beta$.
\begin{eqnarray*}
\langle [\jjinv u_\beta,X_\beta],\tilde{H}_\beta\rangle &=&\langle j(\tilde{H}_\beta)\jjinv u_\beta, X_\beta\rangle\\
&=&\langle \jjinv j(\tilde{H}_\beta)u_\beta,X_\beta\rangle \mbox{ since }Z,\tilde{H}_\beta\in\carto\\
&=& -\langle j(\tilde{H}_\beta)u_\beta,\jjinv X_\beta\rangle\\
&=&\langle u_\beta, j(\tilde{H}_\beta)\jjinv X_\beta\rangle\\
&=&\langle j(\tilde{H}_\beta)\jjinv X_\beta, u_\beta\rangle\\
&=&\langle [\jjinv X_\beta,u_\beta], \tilde{H}_\beta\rangle
\end{eqnarray*}

Thus $[\jjinv X_\beta, u_\beta]=[\jjinv u_\beta, X_\beta]$ for all $\beta\in \Phi\subseteq \Lambda$.\qed
\end{proof}

We now complete the proof of Proposition \ref{maxrank}.
For the 2 cases for $\xi(s)$, using the commutation relations of Section \ref{chap4}, we find:
\begin{eqnarray*} 
dF_Z(u_0)&=&\omega \left(u_0+\frac{1}{\alpha}[u_0,\jjinv X_1]\right)\in\ker j(Z)\oplus \sum_{\beta\in\Lambda}\g_0(\beta)\\
dF_Z(u_\beta)&=&\frac{\omega}{\alpha} \left([X_0,\jjinv u_\beta]+[\jjinv u_\beta, X_\beta]\right)\in\g_0(\beta)\oplus \langle H_\beta\rangle
\end{eqnarray*}

By hypothesis, $\dim_{\mathbb{C}}V_\beta\geq 2$ and hence $\dim_{\mathbb{R}}U_\beta\geq 4$ for all $\beta\in\Phi\subseteq\Lambda$ by (4) of Proposition \ref{eb4.3}.  Additionally, $\dim \g_0(\beta)=2$ by Proposition \ref{eb4.5}.  
By Remark \ref{remark5.5.1} of Section \ref{sec5.4}  there exists an element $u_\beta\in U_\beta$ such that $\mbox{ad}(X_0)(\jjinv u_\beta)=0$ and $\mbox{ad}(X_\beta)(\jjinv u_\beta)\neq 0$.  Hence $dF_Z(u_\beta)=\frac{\omega}{\alpha}[\jjinv u_\beta,X_\beta]$ is a nonzero element of $\langle H_\beta \rangle$ and
since $\langle H_\beta \rangle$ is 1-dimensional, it follows that $\langle H_\beta\rangle\subseteq\imdf$.  Since this is true for each $\beta\in\Lambda$,
we find that $\carto=\bigoplus_{\beta\in\Lambda}\langle H_\beta\rangle \subseteq \imdf$.

Now let $u_\beta\in U_\beta$ be arbitrary.  Since $\langle H_\beta\rangle \subseteq \imdf$ it follows from the expression for $dF_Z(u_\beta)$ that $\frac{\omega}{\alpha}[X_0,\jjinv u_\beta]\in\imdf$.  Hence $\imdf\supseteq \mbox{ad}X_0(U_\beta)=\goo(\beta)$ for $\beta\in\Lambda$ by (1) of Remark \ref{remark5.5.1}.  Therefore $\displaystyle\goo=\carto\oplus \sum_{\beta\in\Lambda}\goo(\beta)\subseteq \imdf$.  In particular, $\displaystyle \frac{1}{\alpha}[u_0,\jjinv X_1]\in\sum_{\beta\in\Lambda}\goo(\beta)\subseteq \imdf$ for all $u_0\in U_0$ and $X_1\in U_1$.  Then from the expression for $dF_Z(u_0)$ we see that $u_0\in\imdf$ for all $u_0\in U_0$.  Note that $U_0=\ker j(\alpha Z)=\ker j(Z)$ since $Z$ is super regular.  Thus $\w_Z=\mbox{span}\{\ker j(\alpha Z),\z\}=$ \\ $\mbox{span}\{U_0,\goo\}\subseteq \imdf$ allowing us to conclude that the first hit map has maximal rank at every point of $\n_Z^*$.\qed
\end{proof}

\section{Main Result}
\label{chap6}

\subsection{Admissible $\goo-$modules}
\label{sec6.1}

We consider the metric 2-step nilpotent Lie algebra $\{\n=U\oplus\goo, \ip\}$ where $\goo$ is a compact real semisimple  Lie algebra and $U$ is a finite dimensional $\goo-$module.  In Section \ref{chap4} we described the root space decomposition of $\goo$ using the corresponding complex Lie algebra $\g=\goc$, finite dimensional $\g-$module $V=U^\mathbb{C}$ and the weight space decomposition of $U$.  Recall that for the first hit map to have maximal rank, the roots of $\g$ must be weights of the $\g-$module $V$ of complex multiplicity at least two.  The list of highest weights $\lambda$ such that this condition does not hold for $V(\lambda)$ is $L'$.  We found this list by reducing to the case $\g$ simple and $V$ irreducible.  We now explain this reduction and use the list  $L'$ to define an admissible $\goo-$module $U$ for our main result.

For any compact real semisimple Lie algebra $\goo$, we write $\goo$ as a direct sum of simple ideals, $\goo=\g_1\oplus\cdots\oplus\g_m$.  If $U$ is a finite dimensional $\goo-$module, then we consider $U$ as a direct sum of irreducible $\g_i-$modules, $U=U_{i,1}\oplus\cdots\oplus U_{i,n_i}$ for each $i=1,\dots,m$.  For each simple ideal $\g_i$ and each irreducible real submodule $U_{i,j}$ we find an associated highest weight $\lambda_{i,j}$ as described below.  
\begin{definition}We say that $U$ is an \emph{admissible $\g_i$-module} if the family of highest weights $\{\lambda_{i,j}\}_{j=1}^{n_i}$ is not contained in $L'$.  We define $U$ to be an \emph{admissible $\goo-$module} if $U$ is an admissible $\g_i-$module for each simple ideal $\g_i$.  
\end{definition}
Recall that $\ker(\goo)=\{0\}$, where $\ker(\goo)=\{X\in U| ZX=0 \mbox{ for all } Z\in\goo\}$.
\paragraph{Admissibility and the root multiplicity condition}
Let $\Phi$ be the set of roots of $\g=\goc$ and $\Lambda$ the set of weights of $\g-$module $V=U^\mathbb{C}$.  To use the commutation relations for $\n=U\oplus\goo$ found in \ref{sec4.6}, all roots of $\g$ must be weights of $V$, i.e. $\Phi\subseteq \Lambda$.  In addition, to prove the maximal rank of the first hit map, these weights that  are also roots must have complex multiplicity at least two.  We prove that if $U$ is an admissible $\goo-$module as defined above then the roots of $\g=\goc$ are weights of $V=U^\mathbb{C}$ of multiplicity greater than or equal to two.  

\begin{proposition} If $U$ is an admissible $\goo-$module, then all roots of $\g=\goc$ are weights of the $\g-$module $V=U^\mathbb{C}$ of multiplicity at least two.\label{admiss}\end{proposition}

\begin{proof}
As above, let $\goo=\g_1\oplus\cdots\oplus\g_m$ where each $\g_i$ is a simple real Lie algebra.  Then $\g=\goc=\g_1^\mathbb{C}\oplus\cdots\oplus\g_m^\mathbb{C}$ where $\g_i^\mathbb{C}$ is a simple complex Lie algebra for each $i$.  Let $\Phi$ be the set of roots of $\g$.  By Corollary \ref{cor14.1} $\Phi=\Phi_1\cup\cdots\cup\Phi_m$ where $\Phi_i$ is an irreducible set of roots of $\g_i^\mathbb{C}$ for each $i$.  Thus it suffices to show that the roots $\Phi_i$ of $\g_i^\mathbb{C}$ are weights of $V$ of multiplicity at least two for each $i=1,\dots,m$.

Considering the $\goo-$module $U$ as a direct sum of irreducible $\g_i-$modules for each $i$, $U=U_{i,1}\oplus\cdots\oplus U_{i,n_i}$, we  let $V_{i,j}=U_{i,j}^\mathbb{C}$ be a $\g_i^\mathbb{C}-$module for each $i,j$.  Then from Proposition \ref{recx} we have two possibilities.  Either $V_{i,j}$ is an irreducible $\g_i^\mathbb{C}-$module or $V_{i,j}=W_{i,j}+J(W_{i,j})$ where $W_{i,j}$ is an irreducible $\g_i^\mathbb{C}-$module.  We consider each case separately.

In the first case let $\lambda_{i,j}$ be the highest weight of the $\g_i^\mathbb{C}-$module $V_{i,j}$.  Since $U$ is an admissible $\g_i-$module, there is at least one highest weight $\lambda_{i,j}$ such that $\lambda_{i,j}\notin L'$, for some $1\leq j\leq n$.  Thus the roots of $\g_i^\mathbb{C}$ are weights of $V_{i,j}$ of multiplicity at least two.

In the second case, let $\lambda_{i,j}$ be the highest weight of the $\g_i^\mathbb{C}-$module $W_{i,j}$.  Since $U$ is an admissible $\g_i-$module, there is a highest weight $\lambda_{i,j}$ such that $\lambda_{i,j}\notin L'$ for some $1\leq j\leq n$.  Thus the roots of $\g_i^\mathbb{C}$ are weights of $W_{i,j}$ of multiplicity at least two.  Since $V_{i,j}=W_{i,j}+J(W_{i,j})$ by properties of weight spaces, the roots are also weights of $V_{i,j}$ of multiplicity at least two.

In both cases the roots of $\g_i^\mathbb{C}$ are weights of $V_{i,j}$ of multiplicity at least two for some $j$, thus the roots $\Phi_i$ of $\g_i^\mathbb{C}$ are weights of $V=V_{i,1}\oplus\cdots\oplus V_{i,n_i}$ of multiplicity at least two.  Also since $U$ is an admissible $\goo-$module, this is true for each $\g_i$, $i=1,\dots,m$.  Thus the roots $\Phi=\Phi_1\cup\cdots\cup\Phi_m$  of $\g$ will be weights of $V$ of multiplicity at least two.\qed
\end{proof}

\subsection{The Main Result}
\label{sec6.2} 

Let $\goo$ be a compact, semisimple Lie algebra  and $U$ an admissible $\goo-$module.  We consider the context in which $\n$ is the metric 2-step nilpotent Lie algebra $\n=U\oplus\g_0$ defined in \ref{sec4.1} with associated simply connected Lie group $N$.  We will now prove our main result, that for any admissible $\goo-$module, the nilmanifold $\ngam$ has the density of closed geodesics property for any lattice $\Gamma$ associated to a Chevalley rational structure of $\g=\goc$.

\begin{theorem}  Let $U$ be an admissible $\g_0$-module for a compact semisimple  Lie algebra $\goo$.  Let 
$\beta$ be a rational basis of $\n=U\oplus \g_0$ determined by a Chevalley basis of $\g=\g_0^\mathbb{C}$ as in \ref{sec4.3}.
Then $\Gamma\backslash N$ satisfies DCG for every lattice $\Gamma$ in $N$ determined by $\beta$.
\label{mythm}
\end{theorem}
\paragraph{Generic vectors}
\label{sec6.3}
An element $\xi=X+Z\in U\oplus \goo=\n$ is said to be generic if $X$ has a nonzero 
 component in each eigenspace of $j(Z)^2$, including $\ker j(Z)$.  Note that $\ker j(Z)\neq \{0\}$ for every $Z\in\goo$ since the weight space $U_0$ is assumed to be nonzero by \ref{sec3.3} and Proposition \ref{eb4.3}.  It is not difficult to see that the generic vectors of $\n$ form a dense open subset of $\n$ (cf. \cite[Proposition 1.19]{gm}).

\paragraph{Reduction to geodesics starting at the identity}
\label{sec6.4}
The following result allows us to reduce our study of geodesics in $N$ to those geodesics $\gamma(t)\in N$ such that $\gamma(0)=e$.

\begin{lemma}\cite[Lemma 13]{demeyer}   Let $\Gamma$ be a lattice in $N$.  Suppose a dense set of geodesics starting
at the identity of $N$ project to closed geodesics in $\Gamma \backslash N$.  Then for a dense
set of points $n\in N$, the set of geodesics with $\gamma(0)=n$ which project to closed
geodesics in $\Gamma \backslash N$ are dense in the set of all geodesics starting at $n$.\label{demeyer13}\end{lemma}
\paragraph{The $m^{th}$ hit map $F_Z^m$}
\label{sec6.5}
Let $Z\in\z_{res}$, the set of resonant vectors in $\z=\goo$, the center of $\n=U\oplus\goo$. The sets $\n_Z$ and $\w_Z$ remain as defined above.

Let $\xi=X+\alpha Z\in\n_Z$.  Since $Z$ is resonant, there exists a smallest value $\omega>0$ such that $e^{\omega j(\alpha Z)}=Id$ on $U$.  Let $\gamma_\xi(t)$ denote the geodesic in $N$ with $\gamma_\xi(0)=e$ and $\gamma_\xi'(0)=\xi\in\n=T_eN$.  Recall from \ref{sec5.1} that the first hit map is defined by $F_Z(\xi)=\log(\gamma_\xi(\omega))$.

For any integer $m\geq 1$ we define the $m^{th}$ hit map  $F_Z^m:\n_Z\rightarrow \w_Z$ by $F_Z^m(\xi)=\log(\gamma_\xi(m\omega))$ where $\omega>0$ is as determined above.
It follows immediately from the formula for $F_Z$ above and the derivation of this formula that $F_Z^m(\xi)=mF_Z(\xi)$ for any $\xi\in\n_Z$ and $m\geq 1$.  In particular $F_Z^m(A)=mF_Z(A)$ for any subset $A\subseteq \n_Z$.

By Remark \ref{remark4.1} in Section \ref{sec5.1},  $\w_Z$ is a rational subalgebra of $\n$ if $Z\in\z$ is rational and resonant.  Hence if $W_Z=\exp(\w_Z)$, then $\Gamma\cap W_Z$ is a lattice in $N$ for a rational resonant vector $Z$ by Proposition \ref{cglattice}. By Remark \ref{remark1.3.4} of Section \ref{sec1.2} any open ball in $\n$ of large enough radius intersects $\log(\Gamma\cap W_Z)$. 

If $Z\in\z$ is rational, resonant and super regular, then $F_Z$ has maximal rank on a dense open subset $\n_Z^*$ of $\n_Z$ by Corollary \ref{maxrankg}.  In particular, if $A$ is a nonempty open subset of $\n_Z$, then $F_Z(A)$ contains an open subset $U$ of $\w_Z$.  It follows from the remark above that if $m\in\mathbb{Z}^+$ is sufficiently large, then $F_Z^m(A)=mF_Z(A)$ contains $mU$ and hence a nonzero element of the lattice $\log(\Gamma\cap W_Z)$.

\label{sec6.6}

\begin{proof}[Proof of Theorem \ref{mythm}]

Let\\ $\n'=\{\xi=X+Z\in\n| \xi\neq 0, \xi \mbox{ is generic and } Z \mbox{ is super regular}\}$.  It follows that $\n'$ is a dense open subset of $\n$ since it is the intersection of the dense open subsets of nonzero, generic and super regular vectors.

To prove the DCG result for $\ngam$, first note that by Lemma \ref{demeyer13}, it suffices to consider geodesics $\gamma_\xi(t)\in N$ with $\xi\in\n'\subseteq\n\cong T_eN$.  We show that for any open subset $\mathcal{O}$ of $\n'$, there exists an element $\xi^*\in\mathcal{O}$ such that $\gamma_{\xi^*}(t)$ projects to a closed geodesic in $\ngam$. To show this, it is enough to show that there exist $\phi\in\Gamma$, $\xi^*\in\mathcal{O}$ and $\omega^*>0$ such that $\phi\cdot\gamma_{\xi^*}(t)=\gamma_{\xi^*}(t+\omega^*)$ for all $t\in\mathbb{R}$ by the discussion following Lemma \ref{onefourtwo}.

Let $\mathcal{O}$ be any open subset of $\n'$, and let $\g_\mathbb{Q}$ be the dense subset of $\goo$ defined in \ref{sec4.4} that consists of resonant, rational vectors of $\goo$.  If $\mathcal{O}_\mathbb{Q}=\{X+Z\in\mathcal{O}|Z\in\g_\mathbb{Q}\}$, then $\mathcal{O}_\mathbb{Q}$ is dense in $\mathcal{O}$.  

Now let $\xi=X+Z\in\mathcal{O}_\mathbb{Q}\subseteq \n'$.  
Since $Z$ is a resonant vector, we can define a first hit map $F_Z:\n_Z\rightarrow \w_Z$.  By Corollary \ref{maxrankg} $F_Z$ has maximal rank on a dense open subset $\n_Z^*$ of $\n_Z$ since $Z$ is a super regular element of $\g_\mathbb{Q}$.  The set $A=\mathcal{O}\cap\n_Z$ is clearly an open subset of $\n_Z$.  The set $A$ is nonempty since it contains $\xi$; $\xi$ is generic since it lies in $\mathcal{O}\subseteq \n'$ and hence $\xi$ also lies in $\n_Z$ by the definition of $\n_Z$.  From the discussion above of the $m^{th}$ hit map it follows that if $m\in\mathbb{Z}^+$ is sufficiently large, then there exists a nonidentity lattice point $\phi\in\Gamma\cap W$ such that $F_Z^m(A)=mF_Z(A)$ contains $\log \phi$.  Let $\xi^*=X^*+\alpha Z\in A$ be an element such that $\phi=\exp(F_Z^m(\xi^*))=\gamma_{\xi^*}(m\omega)$ where $\omega>0$ is the first value such that $e^{\omega j(\alpha Z)}=Id$ on $U$.  By Proposition \ref{eb14.3} it follows that $\phi\cdot\gamma_{\xi^*}(t)=\gamma_{\xi^*}(t+m\omega)$ for all $t\in\mathbb{R}$.  Since $\phi\in\Gamma\cap W\subseteq \Gamma$ and $\xi^*\in A=\mathcal{O}\cap \n_Z\subseteq \mathcal{O}\subseteq\n'$ the proof of the theorem is complete by the discussion above.
\end{proof}

\section{Theory of Weights}
\label{chap3}

Here we include results used earlier.

\subsection{Roots are Weights}
\label{sec3.4}

Let $V$ be a finite dimensional irreducible complex $\g-$module. 
We consider the classification of complex simple Lie algebras into types $A_n,$ $B_n,$ $C_n,$ $D_n,$ $E_6,$ $ E_7,$ $E_8,$  $F_4,$ and
$G_2$.  We show that for each of these Lie algebras, almost all dominant weights of $V$
 are highest weights of irreducible representations for which all roots are weights.  Those exceptional cases for which roots are not weights are a subset of  the list $L'$.  We mentioned the following result of \cite{humph} above.

\begin{lemma}\cite[Lemma 10.4.C]{humph}   Let the set of all roots $\Phi$ be irreducible. 
Then at most 2 root lengths occur in $\Phi$, and all roots of a given length are conjugate under 
the action of $\weyl$, the Weyl group.
\label{rootsconj}
\end{lemma}

Thus we can divide our consideration of the complex simple Lie algebras into two cases: those with all roots the same length and those with two root lengths.  The Lie algebras of type $A_n$,
 $D_n$, $E_6$, $E_7$, and $E_8$ have only one root length, all others have both short and long roots.

In the remainder of this paper we let $\alpha_i$ be the simple roots with some fixed ordering and let $\omega_i$ be the fundamental dominant weights.

\begin{proposition} Let $V$ be an irreducible $\g-$module with nontrivial zero weight space $V_0$.  Let $\lambda\in\Lambda(V)^+$ be the highest weight.  Then
\begin{enumerate}
\item $\displaystyle\lambda=\sum_{i=1}^n p_i\alpha_i$ for suitable positive integers $p_i$.
\item If $\mu\in\Lambda(V)$, then $\displaystyle\mu=\sum_{i=1}^n m_i\alpha_i$, $\mi\in\mathbb{Z}$.  Furthermore, if $\mu\in\Lambda^+(V)$, then the integers $\{m_i\}$ are all positive.
\item At least one root of $\g$ is a weight.
\end{enumerate}
\label{atleast1} \end{proposition}

\begin{remark} \label{remark3.4} The result above, especially the first and second statements, plays a key role in the discussion of this section.  For (3) we show in Proposition \ref{allcases} that in all cases each short root is a weight.
In the proof of Proposition \ref{atleast1} we use the following lemma. \end{remark}

\begin{lemma}  Let $\mu\in\Lambda^+$ and suppose that $\displaystyle\mu=\sum_{i=1}^nm_i\alpha_i$ for integers $\{m_i\}$.  Then $m_k>0$ for all $k$.
\label{posint}
\end{lemma}

\begin{proof} 
Since $\mu\in\Lambda^+$ there exist nonnegative integers $\{r_1,\dots,r_n\}$ such that $\displaystyle \mu=\sum_{i=1}^n r_i\omega_i$.  Let $C^{ij}$ denote the inverse of the Cartan matrix $C_{ij}=\langle \alpha_i,\alpha_j\rangle$.  Recall that $\displaystyle\alpha_i=\sum_{j=1}^n C_{ij}\omega_j$ (cf. \cite{humph}), which implies that $\displaystyle\omega_i=\sum_{k=1}^n C^{ik}\alpha_k$.  Hence $\displaystyle\mu=\sum_{i=1}^nr_i\omega_i=\sum_{i=1}^n m_i\alpha_i$, where $\displaystyle m_k=\sum_{i=1}^nr_iC^{ik}$.  From a case by case consideration of the inverse Cartan matrices (see \cite{humph} p 69), one can see that the elements $C^{ik}$ are always positive.  Hence the integers $\{m_1,\dots,m_n\}$ are positive.
\qed
\end{proof}

\begin{proof}[Proof of Proposition \ref{atleast1}]
Let $\lambda\in\Lambda^+(V)$ be the highest weight of $V$.  Then every weight $\mu\in\Lambda(V)$ has the form $\displaystyle\mu=\lambda-\sum_{i=1}^nq_i\alpha_i$, where $\{q_1,\dots,q_n\}$ are suitable nonnegative integers by \cite{humph}.  
\begin{enumerate}
\item Since $\mu=0\in\Lambda(V)$, we obtain $\displaystyle\lambda=\sum_{i=1}^n p_i\alpha_i$, where $\{p_1,\dots,p_n\}$ are nonnegative integers.  In fact, the integers $\{p_1,\dots,p_n\}$ are positive by Lemma \ref{posint}.
\item By (1) and the discussion above, for any $\mu\in\Lambda(V)$ we have $\displaystyle\mu=\lambda-\sum_{i=1}^n q_i\alpha_i=\sum_{i=1}^n(p_i-q_i)\alpha_i$, where $\{p_i,q_i\}$ are integers.  The remainder of (2) follows from Lemma \ref{posint}.
\item Recall that $\displaystyle \g=\cart\oplus \sum_{\beta\in\Phi}\g_\beta$ and 
that $\displaystyle V=V_0\oplus\sum_{\linl}V_\lambda$.  Note that $\cart(V_0)\equiv 0$ since $V_0$ is the zero weight space.  
We also have $\g_\alpha(V_0)\subseteq V_\alpha$ for all $\alpha\in\Phi$ by Lemma \ref{humph107}.  If $V_\alpha=\{0\}$
for all $\alpha$, then $V_0$ is a proper $\g-$submodule of $V$ on which $\g(V_0)\equiv \{0\}$.  However,
we assumed that $V$ was irreducible, so therefore there must be some $\alpha\in\Phi$ for which
$\g_\alpha(V_0)\subseteq V_\alpha\neq \{0\}$.  Thus the root $\alpha$ is also a weight.
\end{enumerate}
\qed
\end{proof}

\subsection{The simple Lie algebra types}
\label{sec3.5}

\begin{proposition} Let $V$ be an irreducible $\g-$module where $\g$ is of type $A_n$, $D_n$, $E_6$, $E_7$,
or $E_8$.  Let zero be a nontrivial weight of $V$.  Then all roots of $\g$ are weights of $V$.
\end{proposition}

\begin{proof}
By Proposition \ref{atleast1}, in each case we have that one root is a weight.  By Lemma \ref{rootsconj} we know that all roots of the same length are conjugate and therefore all of the remaining roots are conjugate to a weight.  Finally we use the fact that the set of weights is invariant under conjugation by the Weyl group to conclude that all roots are weights.  \qed
\end{proof}

Next, we consider the Lie algebras with two root lengths:  $B_n$, $C_n$, $F_4$, and $G_2$.  Let
$\Delta=\{\alpha_1,\dots,\alpha_n\}$ be the base of simple roots of the Lie algebra of rank $n$.

Let $\mu_1$ and $\mu_2$ denote the highest short and long roots respectively as listed in Table \ref{domwts}. Where there are two root lengths, the highest short root $\mu_1$ is listed first. See also Section 12.2 of \cite{humph}.  The roots $\mu_1$ and $\mu_2$ lie in $\Lambda^+$ by Section 13.2 of \cite{humph}.  By Lemma \ref{rootsconj} $\weyl(\mu_1)$  is the set of all short roots and $\weyl(\mu_2)$ is the set of all long roots.  Hence $\Phi=\weyl(\mu_1)\cup\weyl(\mu_2)$.

\begin{table}[h]
\caption{Highest short and long roots}
\centering
\label{domwts}
\begin{tabular}{|c|c|} \hline
Classical Lie algebra type &  \\ \hline
$A_n$ & $\alpha_1+\alpha_2+\cdots+\alpha_n$\\ \hline
$B_n$ & $\alpha_1+\alpha_2+\cdots +\alpha_n$ \\ 
 	& $\alpha_1+2\alpha_2\cdots +2\alpha_n$ \\ \hline
$C_n$ & $\alpha_1+2\alpha_2+2\alpha_3+\cdots+2\alpha_{n-1}+\alpha_n$\\
	& $2\alpha_1+2\alpha_2+2\alpha_3+\cdots+2\alpha_{n-1}+\alpha_n$\\ \hline
$D_n$ & $\alpha_1+2\alpha_2+\cdots +2\alpha_{n-2}+\alpha_{n-1}+\alpha_n$\\ \hline
$E_6$	& $\alpha_1+2\alpha_2+2\alpha_3+3\alpha_4+2\alpha_5+\alpha_6$\\ \hline
$E_7$	& $2\alpha_1+2\alpha_2+3\alpha_3+4\alpha_4+3\alpha_5+2\alpha_6+\alpha_7$\\ \hline
$E_8$	& $2\alpha_1+3\alpha_2+4\alpha_3+6\alpha_4+5\alpha_5+4\alpha_6+3\alpha_7+2\alpha_8$\\ \hline
$F_4$ & $\alpha_1+2\alpha_2+3\alpha_3+2\alpha_4$ \\ 
	&   $2\alpha_1+3\alpha_2+4\alpha_3+2\alpha_4$ \\ \hline
$G_2$ & $2\alpha_1+\alpha_2$ \\ 
	& $3\alpha_1+2\alpha_2$ \\ \hline
\end{tabular}
\end{table}

\begin{proposition}  Let $\mu_1$ and $\mu_2$ be the highest short root and the highest long root respectively. Let $\linl^+(V)$ be the highest weight.  Then 
\begin{enumerate}
\item $\Lambda(V)$ contains the set of short roots $\weyl(\mu_1)$.
\item If $\lambda \succeq \mu_2$, then $\Lambda(V)$ contains all roots $\Phi$.
\end{enumerate}
 \label{allcases} \end{proposition}

\begin{proof} %
\begin{enumerate}
\item It suffices to show that $\mu_1\in\Phi\cap \Lambda(V)$ since both $\Phi$ and $\Lambda(V)$ are invariant under $\weyl$.  By (3) of Proposition \ref{atleast1} there exists $\alpha\in\Phi\cap\Lambda(V)$, and $\weyl(\alpha)$ contains either $\mu_1$ or $\mu_2$ since $\Phi=\weyl(\mu_1)\cup\weyl(\mu_2)$.  If $\mu_1\in\weyl(\alpha)\subseteq\Phi\cap\Lambda(V)$, then we are done.  Note that $\mu_2\succeq\mu_1$ by inspection of Table \ref{domwts}.  If $\mu_2\in\weyl(\alpha)\subseteq\Phi\cap\Lambda(V)$, then $\mu_1\in\Phi\cap\Lambda(V)$ by Lemma \ref{humphB} since $\Lambda(V)$ is a saturated set of weights.
\item If $\lambda\succeq \mu_2$, then since $\mu_2\succeq\mu_1$ it follows from Lemma \ref{humphB} that $\Lambda(V)$ contains both $\mu_1$ and $\mu_2$.  Hence $\Lambda(V)$ contains $\weyl(\mu_1)\cup\weyl(\mu_2)=\Phi$.
\end{enumerate}\qed
\end{proof}

\begin{proposition} Let $V$ be an irreducible $\g-$module where $\g$ is of type $B_n$, $F_4$, and $G_2$, and let zero be a nontrivial weight of $V$.  If
the highest weight is anything except the following
\begin{center} 
\vspace{0.15in}
\begin{tabular}{|c|c|} \hline
$B_n$ &  $\alpha_1+\cdots +\alpha_n$ \\ \hline
$F_4$ & $\alpha_1+2\alpha_2+3\alpha_3+2\alpha_4$ \\ \hline
$G_2$ & $2\alpha_1+\alpha_2$ \\ \hline
\end{tabular}
\end{center}
then all roots are weights.  In each exceptional case, the short root is a weight, but the
long root is not. \label{roots}
\end{proposition}

Note that these are the highest short root in each case.  We consider the case $C_n$ separately because there are several highest dominant
weights for which not all roots are weights.    We will return to discuss this case after dealing with the simpler ones.

Recall that a weight $\lambda=\mone\alpha_1+\cdots+m_n\alpha_n$ is a dominant weight if and only if $\langle \lambda,\ai\rangle\geq 0$ for each $\ai$ where $(\langle \alpha_i,\alpha_j\rangle)$ is the Cartan matrix.

\begin{proof}[Proof of Proposition \ref{roots}] We consider each class of Lie algebras separately.  We also use repeatedly that if $\lambda\in\Lambda^+(V)$, then $\displaystyle \lambda=\sum_{i=1}^n m_i\alpha_i$, where the integers $\{m_1,\dots,m_n\}$ are positive (Lemma \ref{posint}) and  $\langle \lambda, \alpha_i\rangle \geq 0$.

\begin{enumerate}
\item
\emph{Lie algebras of type $B_n$} 
As above, the weight $\lambda = m_1\alpha_1+\cdots+m_n\alpha_n$, $m_i\in\mathbb{Z}$, $\mi>0$ for $i=1,\dots,n$ is a dominant weight if and only if $\langle \lambda, \alpha_i\rangle \geq 0$.  Equivalently, $\lambda$ is a dominant weight if the following inequalities hold:
\begin{eqnarray}
m_2 &\leq& 2m_1\\
m_{i-1}+m_{i+1} & \leq & 2m_i \ \ \ i=2,\dots n-1\\
m_{n-1} & \leq & m_n.
\end{eqnarray}

\begin{sublemma} If $m_2\geq 2$, then $m_i\geq 2$ for all $i\geq 2$.\label{seven}\end{sublemma}
\begin{proof} By (2) and our hypothesis, $2+m_4\leq m_2 + m_4\leq 2m_3$ and since $m_4\neq 0$, $m_3 \geq 2$.  Next, suppose that $m_i\geq 2$ for all $2\leq i \leq N$ for some number $N< n-1$.  
Then by (2) $2+m_{N+2}\leq m_N+m_{N+2}\leq 2m_{N+1}$.  Since $m_{N+2}>0$, $m_{N+2}\geq 1$.
Thus $m_{N+1}\geq 2$.  Hence $m_i\geq 2$ for $2\leq i\leq n-1$ and $m_n\geq 2$ by (3).\qed\end{proof}

\begin{sublemma} Let $\lambda=m_1\alpha_1+\cdots m_n\alpha_n\in\Lambda^+(V)$.  If $m_1=1$, 
then either $m_i=1$ for $1\leq i\leq n$ or $m_i\geq 2$ for $2\leq i \leq n$.\label{eight}\end{sublemma}
\begin{proof}
Suppose that $m_1=1$.  By inequality (1), $m_2\leq 2m_1=2$.  Then $m_2=1$ or $m_2=2$.  If $m_2=1$, then by (2), $m_3=1$ and by induction on (2), $m_i=1$ for $2\leq i\leq n$.  If $m_2\geq 2$, then by Sublemma \ref{seven}, $m_i\geq 2$ for all $i\geq 2$. \qed\end{proof}

This implies that if $\lambda$ is any dominant weight, then either $\lambda \succeq \mu_2$ or $\lambda=\mu_1$.
Then by Proposition \ref{allcases} and Sublemma \ref{eight}, for a highest
weight $\lambda=m_1\alpha_1+\cdots m_n\alpha_n\neq \mu_1$ all roots are weights, proving Proposition \ref{roots} for the $B_n$ case.

\item \emph{Lie algebras of type $F_4$}
As in the $B_n$ case, we find necessary inequalities for  $\lambda=\mone\aone+\cdots+\mfour\afour$, $\mi>0$ for all $i$, to be a dominant weight:
\setcounter{equation}{0}
\begin{eqnarray}
\mtwo &\leq & 2\mone\\
\mone + \mthree &\leq & 2\mtwo\\
\mfour + 2\mtwo &\leq & 2\mthree\\
\mthree &\leq& 2\mfour
\end{eqnarray}

\begin{sublemma}  If $\lambda=\mone\aone+\cdots+\mfour\afour$ is any dominant weight, then $\mtwo\geq 2$, $\mthree\geq 3$ and $\mfour\geq 2$.\label{frequire}\end{sublemma}

\begin{proof}  It follows from inequality (3) that $\mthree\geq 2$.  From (2) we then see that $\mtwo\geq	2$ also, and hence $m_3\geq 3$ by (3).  Then it follows by (4) that $\mfour\geq 2$.\qed \end{proof}

\begin{sublemma}  If $\mone=1$, then $\lambda=\mu_1$.\label{fmone}\end{sublemma}

\begin{proof}  Let $\mone=1$, then by (1) $\mtwo\leq 2\mone =2$, implying by Sublemma \ref{frequire} that $\mtwo=2$.  Then by (2) $1+\mthree\leq 2\mtwo=4$, hence by Sublemma \ref{frequire} $\mthree=3$.  We conclude by (3) and (4) then that $\mfour=2$.  Thus $\lambda=\mu_1$.\qed\end{proof}

\begin{sublemma} If $\mone=2$, then $\lambda\succeq\mu_2$.\label{fmtwo}\end{sublemma}

\begin{proof}  Let $\mone=2$.  Then by (1) $\mtwo\leq 4$, so $\mtwo=2,3,$ or $4$.  
\begin{enumerate}
\item Let $\mtwo=2$.  Then by (2) and \ref{frequire}, $2+3\leq \mone+\mthree \leq 4$, an obvious contradiction.  Thus if $\mone=2$, then $\mtwo\neq 2$.
\item Let $\mtwo=3$.  By (2) $2+\mthree=\mone+\mthree\leq 2\mtwo=6$, so $\mthree\leq 4$.  By (3) $\mfour+6=\mfour+2\mtwo \leq 2\mthree$, so $\mthree\geq 4$.  Thus $\mthree=4$.  By (3) $\mfour\leq 2$, but  $\mfour\geq 2$ by (4).  Hence $\mfour=2$ and $\lambda=\mu_2$.
\item Suppose $\mtwo=4$. By (3) $\mfour+8=\mfour+2\mtwo\leq 2\mthree$, so clearly $\mthree\geq 5$.  By (4) $5\leq \mthree\leq 2\mfour$ gives $\mfour \geq 3$.  Hence $\lambda\succ \mu_2$.
\end{enumerate}\qed
\end{proof}

This proves Proposition \ref{roots} for a Lie algebra of type $F_4$.  The sublemmas \ref{fmone} and \ref{fmtwo} together show that any dominant weight $\lambda$ has the property that $\lambda=\mu_1$ or $\lambda\succeq \mu_2$.  Therefore by (2) of Proposition \ref{allcases} all roots will be weights for any dominant weight $\lambda\succeq \mu_2$.

\item \emph{Lie algebras of type $G_2$}
  The weight $\lambda=\mone\aone+\mtwo\atwo$ is a dominant weight if and only if the following inequalities are satisfied:
\setcounter{equation}{0}
\begin{eqnarray}
3\mtwo&\leq&2\mone\\
\mone&\leq&2\mtwo
\end{eqnarray}

Note that $\mone\geq 2$ by (1) since $\mtwo\geq 1$.
Suppose $\mone=2$, then by (1) $3\mtwo\leq 4$ and by (2) $2\leq 2\mtwo$ giving $\mtwo=1$ and $\lambda=\mu_1$.  This is the only dominant weight with $\mone=2$.
Next, suppose that $\mone=3$.  Then by (1) $3\mtwo\leq 6$ and by (2) $3\leq 2\mtwo$, giving $\mtwo=2$ and $\lambda=\mu_2$.
If $m_1\geq 4$, then $\mtwo\geq 2$ by (2) and it follows that $\lambda\succ \mu_2$.
We conclude then that for $\lambda$ a dominant weight, $\lambda=\mu_1$ or $\lambda\succeq\mu_2$.  By Proposition \ref{allcases} the proof is complete.
\end{enumerate}
\qed
\end{proof}

\begin{center}\textbf{Lie Algebras of type $\mathbf{C_n}$}\end{center}
 In this case the highest short and long roots are $\mu_1=\alpha_1+2\alpha_2+\cdots+2\alpha_{n-1}+\alpha_n$ and 
$\mu_2=2\alpha_1+2\alpha_2+\cdots+2\alpha_{n-1}+\alpha_n$.  We consider two cases for the dominant weight $\lambda=m_1\alpha_1+m_2\alpha_2+\cdots+m_n\alpha_n$:  $m_1\geq 2$ and $m_1=1$.

\begin{proposition} For the Lie algebra $C_n$, if  $\lambda=m_1\alpha_1+m_2\alpha_2+\cdots+m_n\alpha_n$ is a dominant weight then $m_i\geq 2$, $i=2,\dots,n-1$, and $m_n\geq 1$.  If $m_1\geq 2$ then $\lambda\succeq\mu_2$ and all roots are weights.\label{cn}\end{proposition}

\begin{proof}
Let $\lambda = m_1\alpha_1+\cdots+m_n\alpha_n$, $m_i\in\mathbb{Z}$, $m_i>0$, for $i=1,\dots,n$, be
a dominant weight.  Then $\langle \lambda, \alpha_i\rangle \geq 0$ for all $i=1,\dots,n$ if and only if the following inequalities hold:

\setcounter{equation}{0}
\begin{eqnarray}
m_2 &\leq& 2m_1\\
m_{i-1}+m_{i+1} & \leq & 2m_i \ \ \ i=2,\dots n-2\\
m_{n-2}+2m_n & \leq & 2m_{n-1}\\
m_{n-1} & \leq &2 m_n.
\end{eqnarray}

\begin{sublemma}  For a dominant weight $\lambda=\mone\aone+\cdots+\mn\an$, $\mi\geq 2$ for $i=2,\dots,n-1$.  If $\mi=2$ for some $i=3,\dots,n-1$, then $\mi=2$ for all $i=2,\dots,n-1$ and $\mn=1$.\label{cnrequire}\end{sublemma}

\begin{proof}  Since $\mi\neq 0$ for all $i$, we see by inequality (3) that $3\leq m_{n-2}+2\mn\leq 2m_{n-1}$ and therefore $m_{n-1}\geq 2$.  Then by (2) $3\leq m_{n-3}+m_{n-1}\leq 2m_{n-2}$, we have $m_{n-2}\geq 2$.  Induction on (2) yields $m_i\geq 2$ for $i=2,\dots,n-2$.

If $\mi=2$ for some $i=3,\dots,n-2$, then by (2) $m_{i-1}+m_{i+1}\leq 2m_i=4$, but by the first part of the sublemma, $m_{i-1}\geq 2$ and $m_{i+1}\geq 2$, therefore resulting in the equalities $m_{i-1}=m_{i+1}=2$.  Again we repeatedly use (2) to conclude that $m_i=2$ for $i=2,\dots,n-1$.  If $m_{n-1}=2$, we first use (3) to conclude that $m_{n-2}=2$ and then proceed as above.

Lastly, if $m_i=2$ for $i=2,\dots,n-1$, then by (3), $\mn=1$ must hold.\qed
\end{proof}
Clearly then if $m_1\geq 2$ then $\lambda\succeq \mu_2$ by \ref{cnrequire}.  Hence $\Lambda\supseteq \Phi$ by Proposition \ref{allcases}.\qed
\end{proof}

\begin{center}\textbf{Dominant weights of the form $\lambda=\aone+\mtwo\atwo+\cdots+\mn\an$}\end{center}

Next we characterize those dominant highest weights $\lambda$ of $C_n$ such that not all roots are weights.  By Sublemma \ref{cnrequire} these must be  weights of the form  $\lambda=\mone\aone+\cdots+\mn\an$ with $\mone=1$ since otherwise $\lambda\succeq \mu_2$.

\begin{proposition} If $\lambda=\aone+\mtwo\atwo+\cdots+\mn\an$ is a dominant weight, then exactly one of the following must occur:
\begin{enumerate}
\item $m_i=i$ for $1\leq i\leq n-1$, $\mn=\frac{n-1}{2}$ if $n$ is odd and $\mn=\frac{n}{2}$ if $n$ is even.
\item $m_i=i$ for $1\leq i\leq 2N\leq n-1$, $m_j=2N$ for $2N\leq j\leq n-1$ and $m_n=N$.
\end{enumerate}
\label{cnextra}\end{proposition}

\begin{proof}  Follows directly from Lemmas \ref{cnclass}, \ref{cnmore} and \ref{cnagain} below.\qed \end{proof}

\begin{remark}  It is easy to check that cases (1) and (2) of Proposition \ref{cnextra} satisfy the dominant weight inequalities above.  Hence by Proposition \ref{cn} these two cases are the only examples of irreducible  $\g-$modules $V(\lambda)$ such that $m_1=1$ and $\lambda\succeq \mu_1$ but $\lambda \nsucceq \mu_2$.  By Proposition \ref{allcases}, then in these cases $V(\lambda)$ is an irreducible $\g-$module for which not all roots are weights.\label{remarkcn}\end{remark}

\begin{lemma}  Let $\lambda=\aone+\mtwo\atwo+\cdots+\mn\an$ be a dominant weight.  Then $m_k\leq m_{k+1}$ for $1\leq k\leq n-2$.  If equality holds for some $k$, then $m_k=m_j=m_{n-1}=2m_n$ for $k\leq j\leq n-1$.  In particular $m_j$ is even for $k\leq j\leq n-1$. \label{cnclass}\end{lemma}

\begin{proof}
For $k=1,\dots,n-2$ we consider the assertion 
\begin{eqnarray*} (*_k) \ \ m_k &\leq & m_{k+1} \mbox{ and if equality holds then } m_k=m_j=m_{n-1}=2m_n \\
& &  \mbox{ for } k\leq j \leq n-1. \end{eqnarray*}
We prove that $(*_k)$ holds for $1\leq k\leq n-2$ by backward induction on $k$, starting from the case $k=n-2$.

By inequalities (3) and (4) we have $m_{n-2}+m_{n-1}\leq m_{n-2}+2m_n\leq 2m_{n-1}$, which shows that $m_{n-2}\leq m_{n-1}$.  If equality holds, then by (3) and (4) again $m_{n-1}=2m_n$, which proves $(*_k)$ for $k=n-2$.

Suppose now that $(*_k)$ holds for some integer $k$ with $2\leq k \leq n-2$.  We complete the proof by showing that $(*_{k-1})$ holds.  By (2) and $(*_k)$ we have $m_{k-1}+m_k\leq m_{k-1}+m_{k+1}\leq 2m_k$, which implies that $m_{k-1}\leq m_k$.  If equality holds, then $m_k=m_{k+1}$ and by the equality case of $(*_k)$ we obtain $m_{k-1}=m_k=m_j=m_{n-1}=2m_n$ for $k-1\leq j\leq n-1$.  We have proved $(*_{k-1})$. \qed
\end{proof}

\begin{lemma}Let $\lambda=\aone+\mtwo\atwo+\cdots+\mn\an$, and suppose that $\{\mtwo,\dots,m_{n-1}\}$ are all distinct.  Then 
\begin{enumerate}
\item[(a)] $m_k=k$ for $2\leq k\leq n-1$ and
\item[(b)] $2m_n=n-1$ if $n$ is odd and $2m_n=n$ if $n$ is even.
\end{enumerate}
\label{cnmore}
\end{lemma}

\begin{proof}
\begin{enumerate}
\item[(a)]  For an integer $k$ with $2\leq k\leq n-1$ we consider the assertion 
$$(*_k) \ \  m_i= i \mbox{ for } 2\leq i \leq k.$$
We prove by induction on $k$ that $(*_k)$ holds for $2\leq k \leq n-1$.
Since $\mone=1$ it follows that $\mtwo\leq 2$ by inequality (1), but $\mtwo\geq 2$ by Sublemma \ref{cnrequire}.  Hence $\mtwo=2$, which proves $(*_k)$ for $k=2$.
Suppose that $(*_k)$ holds for all integers $k$ with $2\leq k\leq n-2$.  It suffices to prove that $(*_{k+1})$ holds.  From inequality (2) we have $m_{k-1}+m_{k+1}\leq 2m_k$, which implies that $m_{k+1}\leq k+1$ since $m_{k-1}=k-1$ and $m_k=k$ by $(*_k)$.  By Lemma \ref{cnclass} and $(*_k)$ we have $m_{k+1}\geq m_k=k$ and by hypothesis $m_k\neq m_{k+1}$.  We conclude that $m_{k+1}=k+1$, which proves $(*_{k+1})$.
\item[(b)] By part (a) the inequalities (3) and (4) become $n-2+2m_n\leq 2n-2$ and $n-1\leq 2m_n$.  We conclude that $n-1\leq 2m_n\leq n$, and since $m_n\in\mathbb{Z}$, this proves (b).
\end{enumerate}\qed
\end{proof}

\begin{lemma} Let $\lambda=\aone+\mtwo\atwo+\cdots+\mn\an$ and suppose that $\{\mtwo,\dots,m_i\}$ are all distinct but $m_i=m_{i+1}$ for some integer $i$ with $2\leq i\leq n-2$.  Then \begin{enumerate}  
\item[(a)] $m_k=k$ for $2\leq k\leq i$ and
\item [(b)] $i=2N$ for some integer $N$, $m_j=2N$ for $2N\leq j\leq n-1$ and $m_n=N$.
\end{enumerate}
\label{cnagain}
\end{lemma}
\begin{proof}  The proof of (a) is the same as in Lemma \ref{cnmore}, and (b) follows immediately from Lemma \ref{cnclass}. \qed\end{proof}

\subsection{Dimension of Weights} 
\label{sec3.6}

As indicated previously, the list $L'$ in Table \ref{lprime} contains  not only weights for which not all roots are weights, but also those highest weights for which the roots are weights of multiplicity one.  Since each root is conjugate to either the highest short root $\mu_1$ or long root $\mu_2$, each root has the same multiplicity as either $\mu_1$ or $\mu_2$.  We now restrict our discussion to determining for which highest weights $\lambda$ the weights $\mu_1$ and $\mu_2$ are of multiplicity one.  

Let $\g$ be a complex semisimple Lie algebra and let $V=V(\lambda)$ denote a finite dimensional irreducible $\g-$module with highest weight $\lambda$.
Denote the multiplicity of a weight $\mu$ by $K_{\lambda,\mu}$.  The multiplicity of a weight $\mu$ is by definition the dimension of the weight space $V_\mu$ in $V$.  

Since the highest weight $\lambda$ and all its conjugates have multiplicity one, an obvious set of dominant highest weights for which either $\mu_1$ or $\mu_2$ will have multiplicity one, are the weights $\lambda=\mu_1$ or $\lambda=\mu_2$.  Thus Table \ref{domwts} is a subset of the list $L'$.   
Next, using Theorem \ref{bzthm} we will find all other $\g-$modules $V$ for which $K_{\lambda,\mu} =1$.

\begin{definition} For $\g$ simple, define a pair $(\lambda; \mu)$ of weights in $\Lambda^+$ to be \emph{primitive} if $(\lambda-\mu)$ written as the sum of simple roots has all
positive integer coefficients.  
\end{definition}
By \cite{bz} we will be able to reduce to the primitive case to find all weight spaces of dimension one.  Thus the following result will be the basis of our determination of all dominant weights $\lambda$ such that $K_{\lambda,\mu_i}=1$ for $i=1,2$.  
In the notation of \cite{bz}, $\mathbb{Z}_+$ is the set of all nonnegative integers.

\begin{theorem} [\cite{bz}, Theorem 1.3]  All primitive pairs $(\lambda; \mu)$ such that $K_{\lambda,\mu}=1$, up to
isomorphism of Dynkin diagrams, are exhausted by the following list: 
\begin{enumerate}
\item $A_n$ $(n\geq 1)$: $\lambda=l\omega_1$, $\displaystyle \mu=\sum_{1\leq i \leq n}a_i\omega_i$ 
  where $a_i\in\mathbb{Z}_+$ and \\ 
$\displaystyle (l-\sum_{1\leq i\leq n} ia_i)\in (n+1)\mathbb{N}$ 
\item $B_n$ $(n\geq 2)$: $\lambda=l\omega_1$, $\displaystyle \mu=\sum_{1\leq i \leq n}a_i\omega_i$ 
 where $a_i\in\mathbb{Z}_+$
is even and \\
$\displaystyle (l-1)=\sum_{1\leq i \leq n-1}ia_i+\frac{na_n}{2}$
\item $G_2$: $\lambda=l\omega_2$, $\displaystyle \mu=a_1\omega_1+a_2\omega_2$
  where $a_1,\ a_2 \in \mathbb{Z}_+$,
and $3l-1=2a_1+3a_2$
\item $G_2$:   $\lambda= \omega_1$, $\mu=0$.
\end{enumerate}
\label{bzthm}
\end{theorem}

As before we let $\mu_1$ and $\mu_2$ denote the highest short and long roots as determined by a base $\Delta=\{\aone,\dots,\an\}$ of simple positive roots.
 Using the result above, we first find dominant weights $\lambda$ for which $(\lambda; \mu_i)$ is a primitive pair and $K_{\lambda, \mu_i}=1$ for $i=1 \rm{ \  or \ }2$.  Later we consider the case that $(\lambda;\mu_i)$ is not a primitive pair.  
 
\paragraph{Determining primitive weight pairs}
Note that in the case of one root length there is no highest short root, thus we use the notation $\mu=\mu_2$.  We continue with some previous assumptions.  First, recall that the zero weight space of $V$ is assumed to  be nontrivial.  Also, we assume that the highest weights are always dominant weights.  We only need consider those weights for which all roots are weights.

By Theorem \ref{bzthm}, the only Lie algebras that we need to consider when determining the dimension one weight spaces for primitive pairs are $A_n, \ B_n,$ and $G_2$.  We will consider each case in terms of the highest short and long roots $\mu_1$ and $\mu_2$, determining all $\lambda$ such that $(\lambda;\mu_i)$ is a primitive pair and $K_{\lambda,\mu_i}=1$, $i=1,2$.  

\begin{lemma}  Let $\g$ be a Lie algebra of type $A_n$ and let $V$ be an irreducible $\g-$module with highest weight $\lambda$. Let $\mu=\alpha_1+\cdots+\alpha_n$.  Then the only $\lambda\in\Lambda^+$ with $(\lambda;\mu)$ a primitive pair and $\klam=1$ are $\lambda=k(n\alpha_1+(n-1)\alpha_2+\cdots+\alpha_n)$ with $k\in\mathbb{Z}$. \label{ancase}\end{lemma}
\begin{proof}
By Theorem \ref{bzthm}, $\klam=1$ for  $\lambda=l\omega_1$ and $\displaystyle \mu=\sum_{1\leq i \leq n} a_i\omega_i$ for $\displaystyle a_i\in\mathbb{Z}_+$ such that  $(l-\sum ia_i)\in(n+1)\mathbb{N}$.  In the case of $A_n$, the highest long root is $\mu=\alpha_1+\cdots+\alpha_n$.  
 
We recall that the Cartan matrix gives us the relations between the simple roots $\alpha_i$ and the fundamental dominant weights $\omega_i$; i.e. $\alpha_i=\sum_{j=1}^n \langle \alpha_i, \alpha_j\rangle\omega_j$.  Inverting the Cartan matrix describes the fundamental weights $\{\omega_i\}$ in terms of the simple roots $\{\alpha_i\}$.  See (\cite{humph}, p 69) for a precise description.  In the case of $A_n$ we obtain the following:
\begin{eqnarray*}
\alpha_1&=&2\omega_1-\omega_2\\
\alpha_i&=&-\omega_{i-1}+2\omega_i-\omega_{i+1} \mbox{ for }i=2,\dots,n-1\\
\alpha_n&=& -\omega_{n-1}+2\omega_n
\end{eqnarray*}

In particular $\mu=\omega_1+\omega_n$ and in this case $l-\sum ia_i=l-(1+n)\in (n+1)\mathbb{N}$ is equivalent to $l\in(n+1)\mathbb{N}$.  Therefore $\klam=1$ if  $\lambda=k(n+1)\omega_1$ for $k\in\mathbb{Z}$.  Since $\omega_1=\frac{1}{n+1}(n\alpha_1+(n-1)\alpha_2+\cdots +\alpha_n)$, then $\lambda=k(n\alpha_1+(n-1)\alpha_2+\cdots+\alpha_n)$ for $k\in\mathbb{Z}$ belongs to the set $L'$ of inadmissible weights because it does not satisfy the multiplicity condition. \qed
\end{proof}

\begin{lemma}  Let $\g$ be a Lie algebra of type $B_n$ and let $V$ be an irreducible $\g-$module  with highest weight $\lambda$ such that all roots of $\g$ are weights of $V$. Let $\mu_1=\alpha_1+\alpha_2+\cdots+\alpha_n$ and $\mu_2=\alpha_1+2\alpha_2+\cdots+2\alpha_n$.  Then for all primitive pairs $(\lambda; \mu_i)$, $K_{\lambda,\mu_i}\neq 1$ for $i=1,2$.\label{bncase}\end{lemma}

\begin{proof}
By Theorem \ref{bzthm}, $\klam=1$ for  $\lambda=l\omega_1$ and $\displaystyle \mu=\sum_{1\leq i \leq n} a_i\omega_i$ for $\displaystyle a_i\in\mathbb{Z}_+$ even and $\displaystyle l-1=\sum_{1\leq i \leq n-1}ia_i+\frac{na_n}{2}$.   We determine whether there are any dominant weights $\lambda$ satisfying these conditions with respect to $\mu=\mu_1$ or $\mu=\mu_2$.
 
In the case of $B_n$ we have the following relations from the Cartan matrix:
\begin{eqnarray*}
\alpha_1&=&2\omega_1-\omega_2\\
\alpha_i&=&-\omega_{i-1}+2\omega_i-\omega_{i+1} \mbox{ for }i=2,\dots,n-2\\
\alpha_{n-1}&=&-\omega_{n-2}+2\omega_{n-1}-2\omega_n\\
\alpha_n&=& -\omega_{n-1}+2\omega_n
\end{eqnarray*}

Note that  $\mu_1=\omega_1$.

In this case for $\mu_1=\sum_{1\leq i \leq n}a_i\omega_i$,  $a_1=1$ and $a_i=0$ for $i=2,\dots,n$ and therefore the $a_i$ do not satisfy the requirement from Theorem \ref{bzthm} that  all $a_i$ are even.  Thus there are no dominant weights $\lambda$ such that $K_{\lambda,\mu_1}=1$.
Similarly we find that $\mu_2=\omega_2$ does not satisfy the necessary conditions, so we also conclude that there are no dominant weights  $\lambda$ such that $K_{\lambda, \mu_2}=1$.
\qed \end{proof}

\begin{lemma}  Let $\g$ be a Lie algebra of type $G_2$ and let $V$ be an irreducible $\g-$module with highest weight $\lambda$ such that all roots of $\g$ are weights of $V$. Let $\mu_1=2\alpha_1+\alpha_2$ and $\mu_2=3\alpha_1+2\alpha_2$.  The only $\lambda\in \Lambda^+$ with $(\lambda,\mu_i)$ a primitive pair and $K_{\lambda,\mu_i}=1$ is $\lambda=\mu_2$.  For all other primitive pairs $(\lambda; \mu_i)$, $K_{\lambda,\mu_i}\neq 1$ for $i=1,2$.\label{gcase}\end{lemma}

\begin{proof}
The last two cases of Theorem\ref{bzthm} pertain to the $G_2$ case.  Immediately we can conclude that the second one is not relevant to us as we are concerned with the dimension of nonzero weight spaces.  We consider the other case.  The result of Theorem \ref{bzthm} says that $\klam=1$ if $\lambda=l\omega_2$, $\mu=a_1\omega_1+a_2\omega_2$ for $a_1, \ a_2\in \mathbb{Z}_+$ and $3l-1=2a_1+3a_2$.  

We consider this with respect to the relations:
\begin{eqnarray*}
\alpha_1&=&2\omega_1-\omega_2\\
\alpha_2&=&-3\omega_1+2\omega_2\\
\end{eqnarray*}
We observe that $\mu_1=\omega_1$ and $\mu_2=\omega_2$.   Thus in the $\mu_1$ case, the condition $3l-1=2a_1+3a_2$ is equivalent to $3l-1=2$ or $l=1$ and then $\lambda=\omega_2=\mu_2$.  Note that $\mu_2\in L'$ already because it has a dimension one weight space.

In the $\mu_2$ case, the condition $3l-1=2a_1+3a_2$ is equivalent to $3l-1=3$ or $l=4/3$ which is not possible since $\lambda$ has integer coefficients.\qed
\end{proof}

\paragraph{Nonprimitive weight pairs}
\label{sec3.7}
The nonprimitive case occurs when $\zeta$ is a dominant weight such that the difference $\zeta-\mu$ does not contain terms for each simple root; i.e. $\zeta$ and $\mu$ have the same coefficient for at least one simple root.  In this case we must examine each type of complex simple Lie algebra to find all such $\zeta$ such that $(\zeta;\mu)$ is a nonprimitive pair for $\mu$, either the highest short or long root.  Then we will determine if $\kzeta=1$ in each case.  Table \ref{nonprim} contains all such $\zeta$ that fail the dimension requirement.  The proof is contained in \cite{decoste}.  The results are based on the following proposition from \cite{bz}.

\begin{table}[h]
\caption{Nonprimitive Weight Pairs with $K_{\zeta,\mu}=1$}
\centering
\label{nonprim}
\begin{tabular}{|c|c|} \hline
Lie algebra type & dominant highest weight \\ \hline
$A_n$ & $\alpha_1+2\alpha_2$, $2\alpha_1+\alpha_2$ $(n=2)$  \\
	& $\alpha_1+2\alpha_2+\alpha_3$ $(n=3)$  \\ \hline
$B_n$ $(n\geq 2)$ & $\alpha_1+2\alpha_2+m_3\alpha_3$, $m_3\geq 3$, $(n=3)$\\
	& $2\alpha_1+2\alpha_2\cdots +2\alpha_n$ $(n\geq2)$\\ 
	& $\alpha_1+2\alpha_2+3\alpha_3+3\alpha_4+\cdots +3\alpha_n$ $(n\geq3)$\\ \hline
$D_n$ $(n\geq 4)$  &  $\alpha_1+2\alpha_2+2\alpha_3+\alpha_4$, $\alpha_1+2\alpha_2+\alpha_3+2\alpha_4$ $(n=4)$\\ \hline
$G_4$ & $4\alpha_1+2\alpha_2$\\ \hline
\end{tabular}
\end{table}

\begin{proposition} \cite[Proposition 2.4]{bz}  \begin{enumerate}
\item Let $S$ be a subset of simple roots.  Let $\linl$ be an element such that the expansion of the weight $(\lambda-\mu)$ in terms of simple roots involves only elements of $S$.  Then $\klam=K_{p(\lambda),p(\mu)}$.
\item Under the assumptions of part 1, let $S_1,\dots,S_k$ be all the connected components of the set $S$ in the Dynkin diagram of the system of positive roots, and let $\lambda_i=p_{S_i}(\lambda)$ and $\mu_i=p_{S_i}(\mu)$.  Then $\klam=\Pi_{1\leq i\leq k}K_{\lambda_i,\mu_i}$.
\end{enumerate}
\label{bztwo}
\end{proposition}

\subsection{Weight Space Decomposition of a Real $\goo-$module: results from \cite{ebpp}}
\label{sec4.6}

We establish a weight space decomposition of the finite dimensional real $\goo-$module $U$ using the weight space decomposition of the finite dimensional complex $\g-$module $V=U^\mathbb{C}$.  For a more thorough discussion of the material in this section and for proofs of  all of the results, see Section 4 of \cite{ebpp}.

\begin{proposition}\cite[Proposition 6.1]{ebpp} Let $\carto$ be a maximal abelian subalgebra of $\goo$, and let $\cart=\cartoc$ be the corresponding Cartan subalgebra of $\g=\goc$.  Then $i\carto=\cart_{\mathbb{R}}=\{H\in\cart|\alpha(H)\in\mathbb{R} \mbox{ for all } \alpha\in\Phi\}$. \label{eb4.1} \end{proposition}

\begin{proposition}\cite[Proposition 6.2]{ebpp}  If $\linl=\Lambda(V)$, then $-\linl$ and $J(V_\lambda)=V_{-\lambda}$, where $J$ denotes conjugation in $V$ induced by $U$.  In particular, $\dim V_\lambda=\dim V_{-\lambda}$.  Moreover, $J(V_0)=V_0$.\label{eb4.2} \end{proposition}

\begin{proposition}\cite[Proposition 6.3]{ebpp}   For each $\linl$ let $U_\lambda=(V_\lambda\oplus V_{-\lambda})\cap U$.  Let $U_0=V_0\cap U$.  Then
\begin{enumerate}
\item If $H_0\in\carto$, then $H_0(U_\lambda)\subseteq U_\lambda$.
\item $U_\lambda=Re(V_\lambda)=Im(V_\lambda)=Re(V_{-\lambda})=Im(V_{-\lambda})$.  $U_0=Re(V_0)=Im(V_0)$.
\item $U_\lambda^{\mathbb{C}}=V_\lambda\oplus V_{-\lambda}$.  $U_0^{\mathbb{C}}=V_0$.
\item $\dim_{\mathbb{R}}U_\lambda=2\dim_{\mathbb{C}}V_\lambda$. 
\item Let $\Lambda_0=i\Lambda$.  Then 
\begin{enumerate}
\item $\Lambda_0\subseteq \mbox{Hom}(\carto,\mathbb{R})$
\item Let $\linl$ and let $\lambda_0=i\lambda\in\Lambda_0$.  Then $U_\lambda=\{u\in U| H_0^2(u)=-\lambda_0(H_0)^2u \mbox{ for }\\
 H_0\in\carto\}$.
\end{enumerate}
\end{enumerate} \label{eb4.3}
\end{proposition}
\begin{proposition}\cite[Proposition 6.4]{ebpp}   Let $\Lambda^+$ be any subset of $\Lambda$ such that $\Lambda$ is the disjoint union of $\Lambda^+$ and $-\Lambda^+$.  Then
\begin{equation*}
U=U_0+\sum_{\linl^+}U_\lambda \mbox{ \ \ \ (direct sum)}
\end{equation*}
If $\langle \ , \ \rangle$ is any inner product on $U$ such that the elements of $\goo\subseteq \mbox{End}(U)$ are skew symmetric, then the decomposition of $U$ above is an orthogonal direct sum decomposition. \label{eb4.4}
\end{proposition}
Using our root space decomposition for $\g=\goc$, $\displaystyle \g=\cart+\sum_{\alpha\in\Phi}\g_\alpha$, we find the root space decomposition of $\goo$ determined by $\carto$.
\begin{proposition}\cite[Proposition 6.5]{ebpp}   Let $\goo$ be a compact, semisimple real Lie algebra, and let $\carto$ be a maximal abelian subspace of $\goo$.  Then
$$\goo=\carto+\sum_{\alpha\in\Phi}\goo(\alpha)$$ where $\goo(\alpha)=(\g_\alpha\oplus\g_{-\alpha})\cap\goo$, $\g=\goc$ and $\cart=\cartoc$.  Each subspace $\goo(\alpha)$ is 2-dimensional. \label{eb4.5} \end{proposition}
If $\ip=-B_0$, where $B_0$ denotes the Killing form of $\goo$, then the elements of $\mbox{ad}(\goo)$ are skew symmetric on $\goo$ relative to $\ip$. In particular, the decomposition above is orthogonal relative to $\ip$ by Proposition \ref{eb4.4}.
\begin{proposition}\cite[Proposition 6.6]{ebpp}   Let $\displaystyle U=U_0+\sum_{\linl^+}U_\lambda$ (direct sum) be the weight space decomposition of $U$.  Let $\goo(\beta)=(\g_\beta\oplus \g_{-\beta})\cap \goo$ for all $\beta\in\Phi^+$.  Then  $\goo(\beta)(U_\lambda)\subseteq U_{\lambda+\beta}\oplus U_{\lambda-\beta}$ for all $\linl$, $\beta\in\Phi^+$. \label{eb4.6}
\end{proposition}
Typically, we can reduce to the case that $\n=U\oplus\goo$, where $U$ is an irreducible $\goo-$module.  In the next result we relate $U$ and the complex $\g-$module $V=U^{\mathbb{C}}$.    
\begin{proposition}\cite[Proposition 8.1]{ebpp}  Let $U$ be a finite dimensional irreducible real $\goo-$module and let $V=U^\mathbb{C}$.  Then one of the following occurs:
\begin{enumerate}
\item $V$ is an irreducible, complex $\g-$module 
\item If $W$ is any irreducible, complex $\g-$submodule of $V$, then $V=W\oplus J(W)$, where $J:V\rightarrow V$ denotes the conjugation induced by $U$.  Moreover
\begin{enumerate}
\item The map $\phi: W^{\mathbb{R}}\rightarrow U$ given by $\phi(w)=w+Jw$ is a $\goo-$isomorphism between the real $\goo-$modules $W^{\mathbb{R}}$ and $U$.
\item $\Lambda(V)=\Lambda(W)\cup -\Lambda(W)$.
\item $V_\mu=W_\mu\oplus J(W_{-\mu})$ for all $\mu\in\Lambda(V)$ and $V_0=W_0\oplus J(W_0)$.
\end{enumerate}
\end{enumerate}
\label{recx}
\end{proposition}
In case 1 of the above proposition, let $W$ be defined to be $U^\mathbb{C}$, otherwise let it be as defined in case 2.  We call $W$ the irreducible $\g-$module associated to the irreducible $\goo-$module $U$.
\begin{proposition} \cite[Preposition 4.7]{ebpp}  Let $\goo$ be a compact, semisimple real Lie algebra, and let $\carto$ be a maximal abelian subalgebra of $\goo$.  Let $U$ be a real $\goo-$module, and let $\{\n=U\oplus\goo, \langle \ , \ \rangle\}$ be the real 2-step nilpotent Lie algebra defined in Section \ref{sec4.1}. Let $V=U^{\mathbb{C}}$, $\cart=\cartoc$ and $\g=\goc$ and let $\Lambda\subseteq \mbox{Hom}(\cart,\mathbb{C})$ be the set of weights for $V$ determined by $\cart$.  Let $$U=U_0+ \sum_{\linl^+}U_\lambda \ \ \ \mbox{(direct sum)}$$ be the corresponding weight space decomposition of $U$.  Then
\begin{enumerate}
\item $ [U_0,U_\lambda]=\goo(\lambda)\mbox{ for all }\lambda\in\Lambda$
\item $[U_\lambda,U_\mu]= \goo(\lambda+\mu)\oplus\goo(\lambda-\mu)$ for all $\lambda, \mu \in \Lambda$
\item $[U_0,U_0]=\{0\}$
\end{enumerate}
We define $\goo(\lambda)=\{0\}$ if $\lambda\notin\Phi$.
\label{commrel}
\end{proposition}

The notion of a weight vector in the real situation, $\n=U\oplus \goo$ parallels that of the complex case.  We let $\tilde{H}_\lambda\in\carto$ be the unique vector such that $\langle H,\tilde{H}_\lambda\rangle=-i\lambda(H)$ for all $H\in\carto$.  Then $\tilde{H}_\lambda$ is the weight vector determined by $\lambda$.  

\begin{corollary}\cite[Corollary 4.8]{ebpp}   Let $\n=U\oplus\goo$ be as in Proposition \ref{commrel}. Then 
\begin{enumerate}
\item $[U_\lambda,U_\lambda]=\mathbb{R}-\mbox{span}\{\tilde{H}_\lambda\} \mbox{ if } \lambda\in\Lambda, 2\lambda \notin \Phi$
\item$[U_\beta,U_\beta] =\mathbb{R}-\mbox{span}\{\tilde{H}_\beta\} \mbox{ if } \beta\in\Lambda\cap\Phi$ 
\end{enumerate}\label{eb4.8}\end{corollary}



\end{document}